\documentclass[11pt,a4paper,reqno]{article}
\usepackage{amsmath,amssymb,exscale,mathrsfs}
\usepackage[pdftex]{graphicx,color,epsfig}
\usepackage{subfigure}
\usepackage[subfigure]{ccaption}
\usepackage{float}
\usepackage{enumitem}

\usepackage{color}
\definecolor{marin}{rgb}   {0.,   0.3,   0.7} 
\definecolor{rouge}{rgb}   {0.8,   0.,   0.} 
\definecolor{sepia}{rgb}   {0.8,   0.5,   0.} 
\usepackage[colorlinks,citecolor=marin,linkcolor=rouge,
            bookmarksopen,
            bookmarksnumbered
           ]{hyperref}

\newtheorem{lemma}{Lemma}[section]

\newtheorem{theorem}[lemma]{Theorem}
\newtheorem{proposition}[lemma]{Proposition}
\newtheorem{corollary}[lemma]{Corollary}
\newtheorem{remark}[lemma]{Remark}
\newtheorem{example}[lemma]{Example}

\newtheorem{notation}[lemma]{Notation}
\newtheorem{definition}[lemma]{Definition}
\newtheorem{conclusion}[lemma]{Conclusion}
\numberwithin{equation}{section}
\newcommand{\QED}{\mbox{}\hfill \raisebox{-0.2pt}{\rule{5.6pt}{6pt}\rule{0pt}{0pt}} 
          \medskip\par}             
\newenvironment{Proof}{\noindent
    \abovedisplayskip = 0.5\abovedisplayskip
    \belowdisplayskip=\abovedisplayskip{\bfseries Proof. }}{\QED}

\newcommand{\N}{\mathbb{N}}

\newcommand{\R}{\mathbb{R}}

\newcommand{\C}{\mathbb{C}}

\newcommand{\Sp}{\mathbb{S}}

\newcommand{\Z}{\mathbb{Z}}

\newcommand{\diver}{\mbox{\sl div\,}}

\newcommand{\be}{\begin{equation}}
\newcommand{\ee}{\end{equation}}
\newcommand{\bea}{\begin{eqnarray}}
\newcommand{\eea}{\end{eqnarray}}
\newcommand{\bee}{\begin{eqnarray*}}
\newcommand{\eee}{\end{eqnarray*}}
\def\ds{\displaystyle}

\def\eps{\varepsilon}

\def\fref#1{{\rm (\ref{#1})}}

\def\pa{\partial}

\author{Philippe Chartier \thanks{INRIA Rennes, IRMAR and ENS Rennes, IPSO Project Team, 
Campus de Beaulieu, F-35042 Rennes, France. E-mail: Philippe.Chartier@inria.fr} 
  \and 
  Nicolas Crouseilles \thanks{INRIA Rennes, IRMAR and ENS Rennes, IPSO Project Team, 
Campus de Beaulieu, F-35042 Rennes, France. E-mail: Nicolas.Crouseilles@inria.fr} 
 \and 
Mohammed Lemou
 \thanks{CNRS, IRMAR and ENS Rennes, IPSO Project Team, 
Campus de Beaulieu, F-35042 Rennes, France. E-mail: Mohammed.Lemou@univ-rennes1.fr}
}
\title{Averaging of highly-oscillatory transport equations}        
\begin{document}
\maketitle

\begin{abstract}
In this paper, we develop a new strategy aimed at obtaining high-order asymptotic models for transport equations with highly-oscillatory solutions. The technique relies upon recent developments averaging theory for ordinary differential equations, in particular normal form expansions in the vanishing parameter. Noteworthy, the result we state here also allows for the complete recovery of the exact solution from the asymptotic model. This is done by solving a companion transport equation that stems naturally from the change of variables underlying high-order averaging. Eventually, we apply our technique to the Vlasov equation with external electric and magnetic fields. Both constant and non-constant magnetic fields are envisaged, and asymptotic models already documented in the literature and re-derived using  our methodology. In addition, it is shown how to obtain new high-order asymptotic models. 
\end{abstract}

\section{Introduction}
In a large variety of situations, one is confronted to the resolution of a family of transport equations of the form 
\begin{eqnarray} \label{eq:transport}
 \partial_t f(t,y) + F^\eps(y) \cdot \nabla_y f(t,y) = 0, \quad f(0,y)=f_0(y) \in \R,\quad t \in \R, \quad y \in \R^n, 
\end{eqnarray}
indexed by a small positive parameter $\eps$, whose occurrence in {\em real-life} models often lies at the core of numerous theoretical and numerical difficulties encountered in obtaining a(-n) (approximate-) solution. The nature of the difficulties (both theoretical and numerical) triggered by the presence of $\eps$ may vary according to the form of the vector field $y \mapsto F^\eps(y) \in \R^n$. In this article, we shall address the highly-oscillatory situation where it can be  split into two parts
\begin{eqnarray} \label{eq:splitF}
F^\eps(y) = \frac{1}{\eps} \omega (y) \, G(y) + K(y)
\end{eqnarray}
where the flow $(t,y_0) \mapsto \Phi_t(y_0)$ associated with the differential equation 
\begin{eqnarray} \label{eq:eqdiffG}
\dot y(t) = G(y(t)), \quad y(0)=y_0,
\end{eqnarray}
is assumed to be {\bf periodic}, regardless of the specific trajectory (i.e. independently of the initial condition $y_0$ at time $t=0$) and where $y \mapsto \omega(y)$ is a  scalar function {\em bounded from below} by a positive constant. Owing to the $1/\eps$-term in front of the vector field $G$, the solution of the transport equation evolves in a {\em highly-oscillatory regime} as soon as $\eps$ becomes small, which is specifically the regime under investigation here. Since our ultimate goal is the design of {\em high-order  uniformly accurate} numerical methods (i.e. methods whose computational cost and accuracy are not influenced by the value of $\eps$), the identification of the {\rm asymptotic models} is a pre-requisite: this is the task addressed in this work.

Examples of highly-oscillatory equations of the form  (\ref{eq:transport}) are numerous \cite{bostan2,bostan3,bostan5,bostan6,frenod3,frenod1,frenod2,golse}. It is obviously out of the scope of this introductory paper to treat all of them: we will rather concentrate on the following model that will constitute hereafter our {\em target application}, namely the {\bf Vlasov equation with strong magnetic field}
\begin{equation} \label{eq:vlasov}
\partial_t f(t,x,v) + v \cdot \nabla_x f(t,x,v) + \left(E(x) + \frac{1}{\eps} v \times B(x)\right) \cdot \nabla_v f(t,x,v) = 0,
\end{equation}
where $x \in \R^3$ and $v \in \R^3$ denote respectively the spatial and velocity variables, $f: \R \times \R^3 \times \R^3 \mapsto \R$ is the distribution function, i.e. the density of particles at time $t$, position $x$ and velocity $v$, and where $E: \R^3 \mapsto \R^3$ and $B: \R^3 \mapsto \R^3$ are respectively the electric and magnetic fields, assumed to be {\em external} at this stage (i.e. not coupled with $f$ through Maxwell equations for instance). 

Our first objective is this paper is to derive {\em formal} asymptotic models for equation (\ref{eq:transport}) with $F^\eps$ satisfying  (\ref{eq:splitF}) and $\omega \equiv 1$. Rather then  merely obtain the  limit equation where $\eps$ tends to zero, we demand  higher-order terms in powers of $\eps$. 
The methodology we propose relies on recent results from the theory of averaging  for highly-oscillatory ordinary differential equations \cite{perko,sanders}, and more precisely on normal forms obtained as $\eps$-expansions. Such series have been derived with the help of B-series in \cite{chartier4,chartier3,chartier2} or somehow more simply in \cite{murua1,murua2,murua3} with word-series\footnote{Although the effect of truncating the aforementioned formal series has been fully analysed in subsequent papers \cite{chartier3, chartier5}, it is out of the scope of this first paper to present a complete error analysis. This will be the object of a forthcoming paper \cite{chartierfuture}.}. The underlying results we shall lean onto will be presented in Section \ref{sect:averode}, but prior to that, we shall show in Section \ref{sect:split} how the splitting of the vector field $F^\eps$ into two {\em commuting} vector fields naturally leads to {\em two independent transport equations}\footnote{The aim of this section is to introduce the rationale underlying our methodology, i.e. the idea that decomposing the vector field $F^\eps$ in (\ref{eq:transport}) into two commuting vector fields allows to separate the stiff and non-stiff parts of the transport equation.}. The corresponding first result   (for constant $\omega$) will be stated in Section \ref{sect:averkin}.
 
In Section \ref{sect:varyfreq}, we will address the much more involved situation of a varying frequency ($\omega$ non-constant in (\ref{eq:splitF})), which requires to work in an augmented space.  In particular, the main result of this paper will be stated there. It allows to rewrite the original transport equation (\ref{eq:transport}) as a set of four {\em non-stiff} equations for a {\em phase function ($S$)} and a {\em profile function ($h$)}. This procedure is inspired from the recent work \cite{crouseilles}, although the context here is different. The two equations for the profile function are the counterpart of the averaged equation obtained elsewhere in the literature. However, solving the equation for the phase function $S$ allows to recover exactly the complete solution of (\ref{eq:transport}). This part is up to our knowledge completely new. Since we use series-expansions, it is possible to write down explicitly and in a systematic way the terms appearing in the four equations for $S$ and $h$. In   Section \ref{sect:vlasov}, we shall eventually envisage our target application (\ref{eq:vlasov}) and show how to obtain the terms of these developments. Firstly, in Section \ref{sect:constant}, we will consider the case of a constant magnetic field $B(x) \equiv B$ in (\ref{eq:vlasov}) with 
a physical space of dimension two, 
as it appears to be a simple application of the results of Section \ref{sect:averkin}.
Secondly, in Section \ref{sect:vary}, we will address the more involved situation of a varying magnetic field ($B$ non-constant in (\ref{eq:vlasov})), which requires a preliminary treatment of the transport equation, as exposed in Section \ref{sect:varyfreq}. At last, we shall treat equation (\ref{eq:vlasov}) in full generality, i.e. in three dimensions and with a general magnetic field, and compare the equations we obtain with our methodology to results previously published in the literature. 
\section{Decomposition of a transport equation}\label{sect:split}
Let us consider the Liouville equation 
$$
\partial_t f(t,y) + F(y) \cdot \nabla_y f(t,y) = 0,
$$
associated to a split vector field of the form
$$
 F = F_1 + F_2, 
$$
and let us make the fundamental assumption that the Lie bracket of $F_1$ and $F_2$ vanishes, that is to say that 
$$
\forall y \in \R^n, \quad [F_1,F_2](y) := (\partial_y F_1)(y) \;F_2(y) - (\partial_y F_2)(y) \; F_1(y)=0.
$$
This {\em commutation} of vector fields further manifests itself as the commutation of the two flows\footnote{These flows are assumed to be defined for all $t \in \R$ and all $y \in \R^n$ without further notice.} associated with $F_1$ and $F_2$, or as the commutation of the Lie operators associated with $F_1$ and $F_2$. 
More precisely, denoting ${\cal L}_{F_1}$ and ${\cal L}_{F_2}$ the operators defined, for any smooth function $g \in {\cal C}^\infty(\R^n;\R^m)$ by
\begin{align*}
\forall y \in \R^n, \quad {\cal L}_{F_1}(g)(y)  = \partial_y g(y) \, F_1(y) \quad \mbox{ and } \quad \forall y \in \R^n, \quad
{\cal L}_{F_2}(g)(y)  = \partial_y g(y) \, F_2(y),
\end{align*}
we have\footnote{Owing to the general well-known formula ${\cal L}_{F_1} {\cal L}_{F_2} - {\cal L}_{F_2} {\cal L}_{F_1} = {\cal L}_{[F_1,F_2]}$.}
\begin{align} \label{eq:liecom}
{\cal L}_{F_1} {\cal L}_{F_2} = {\cal L}_{F_2} {\cal L}_{F_1},
\end{align}
i.e. more explicitly 
\begin{align*}
\forall g \in {\cal C}^\infty(\R^n;\R^m), \quad {\cal L}_{F_1} \Big({\cal L}_{F_2}(g) \Big)= {\cal L}_{F_2} \Big( {\cal L}_{F_1}(g)\Big).
\end{align*}
%
The method of characteristics immediately gives for any {\em smooth} solution of (\ref{eq:transport}) 
\begin{align}
\forall t  \in \R, \quad f(t,\cdot) = \exp{(-t {\cal L}_{F_1+F_2})} (f_0),
\end{align}
which, owing to relation (\ref{eq:liecom}), can also be written as
\begin{align}
\forall t  \in \R, \quad f(t,\cdot) = \exp{(-t {\cal L}_{F_1})} \exp{(-t {\cal L}_{F_2})} (f_0) = \exp{(-t {\cal L}_{F_2})} \exp{(-t {\cal L}_{F_1})} (f_0).
\end{align}
A somehow {\em natural} step forward now consists in {\em separating the two times} in previous relation and defining the new function with additional variable $\tau$ 
\begin{eqnarray} 
\tilde f(t,\tau,\cdot) = \exp{(-\tau {\cal L}_{F_1})} \exp{(-t {\cal L}_{F_2})} (f_0) = \exp{(-t {\cal L}_{F_2})} \exp{(-\tau {\cal L}_{F_1})} (f_0).
\end{eqnarray}
We are now in position to state the following proposition, 
which shows that the augmented function $\tilde f$ is in fact the unique solution of a system of  two {\bf independent} equations.
\begin{proposition} \label{lem:splittrans}
Consider the  system composed of the following  two transport equations
\begin{eqnarray} \label{eq:transporttau}
\forall (t,\tau,y) \in \R \times \R \times \R^n, \quad  \partial_\tau \tilde f(t,\tau,y) + F_1(y) \cdot \nabla_y \tilde f (t,\tau,y) = 0
\end{eqnarray} 
and 
\begin{eqnarray} \label{eq:transportt}
 \forall (t,\tau,y) \in \R \times \R \times \R^n, \quad  \partial_t \tilde f(t,\tau,y) + F_2(y) \cdot \nabla_y \tilde f (t,\tau,y) = 0,
\end{eqnarray} 
together the with initial condition $\tilde f(0,0,y)=f_0(y)$. If the condition $[F_1,F_2] = 0$
is satisfied, this system has a unique solution, which furthermore satisfies 
$$
\forall (t,y) \in \R \times \R^n, \quad \tilde f(t,t,y)=f(t,y).
$$
\end{proposition}
\begin{Proof}
We first note that, if a solution $\tilde f$ exists, then equations (\ref{eq:transporttau}) and (\ref{eq:transportt}) can be solved in any order. Hence, we can obtain the value of $\tilde f(t,\tau,y)$ by first solving (\ref{eq:transporttau}) for $t=0$ from the initial value $\tilde f(0,0,y)=f_0(y)$ -this furnishes $\tilde f(0,\tau,y)$-  and then by solving (\ref{eq:transportt}) for fixed $\tau$ from this initial value. Insofar as the solution exists, it is thus unique. Now, define 
$$
\tilde f(t,\tau,\cdot) = \exp{(-\tau {\cal L}_{F_1})} \exp{(-t {\cal L}_{F_2})} (f_0) = \exp{(-t {\cal L}_{F_2})} \exp{(-\tau {\cal L}_{F_1})} (f_0).
$$
It is easy to check that it satisfies both (\ref{eq:transporttau}) and (\ref{eq:transportt}) by considering successively the first and the second form. The function $\tilde f$ defined above is thus the unique solution of system (\ref{eq:transporttau}-\ref{eq:transportt}). Finally, 
\begin{align*}
\partial_t \Big(\tilde f(t,t,y) \Big) + F \cdot \nabla_y \tilde f(t,t,y) &= \partial_t \tilde f(t,t,y) +\partial_\tau \tilde f(t,t,y) +F \cdot \nabla_y \tilde f(t,t,y) \\
&= \partial_t \tilde f(t,t,y) + F_1 \cdot \nabla_y \tilde f(t,t,y) + \partial_\tau \tilde f(t,t,y) + F_2 \cdot \nabla_y \tilde f(t,t,y) \\
& = 0.
\end{align*}
The initial condition $\tilde f(0,0,\cdot)=f_0$ and a uniqueness argument then allow to conclude. 
\end{Proof}
%
\section{Averaging of ordinary differential equations in a nutshell} \label{sect:averode}
Since our approach for averaging the transport equation (\ref{eq:transport}) consists in averaging first the characteristics and then rewrite the corresponding Liouville equations, we hereafter recall the main results upon which we shall lean. In this paper, we content ourselves with formal expansions, thus neglecting at this stage the occurrence of error terms. This is justified by the fact that, under appropriate smoothness assumptions, these errors actually become of size $\eps^{n}$ for any fixed $n$, or even exponentially small (i.e. bounded by $Ce^{-C/\eps}$ for some positive constant $C$). A completely rigorous treatment of these error terms for ordinary differential equations can be found for instance in \cite{chartier3}, and their influence in our situation will be analysed in a forthcoming paper \cite{chartierfuture}.
\subsection{A normal form theorem}
Consider the highly-oscillatory differential equation
\begin{eqnarray} \label{eq:ode}
\dot y = F^\eps(y) := \frac{1}{\eps} { G}(y) + { K}(y)
\end{eqnarray} 
i.e. equation (\ref{eq:splitF}) with $\omega \equiv 1$, 
where both vector fields $ G$ and $ K$ are assumed to be smooth\footnote{Either of class $C^k$ or analytic. The precise smoothness assumption determines the type of error bounds, either polynomial or exponential in $\eps$ and is thus not essential here (see \cite{chartierfuture}).}. As already alluded to in the Introduction section,   the fundamental assumption $\bf (H)$ required to go any further is that \\ \\
{\bf 
(H)} $G$ generates a {\bf periodic flow} $\Phi_\tau$, {regardless of the specific trajectory} (i.e. with a period which remains independent of the initial value). By convention, we will suppose here that this period is $2\pi$. \\

Since the Lie bracket of $G$ and $K$ has here no reason to vanish, 
we can not reproduce right away the analysis conducted in previous section. It is precisely the aim of averaging to rewrite $F^\eps$ as the sum of two commuting fields\footnote{At least, this is one way to envisage averaging for ordinary differential equations and this is the point of view adopted both in \cite{chartier2} and in the recent series of papers by Murua and Sanz-Serna \cite{murua1,murua2,murua3}.}. As already emphasized, this is in general possible only up to small error terms, so that the theorem stated below is to be understood in a formal sense. 
\begin{theorem} \label{th:aver}
Suppose that the vector field $F^\eps$ can be split according to equation (\ref{eq:ode}) and that $G$ satisfies assumption $\bf (H)$. Then there exist two vector fields $G^\eps$ and $K^\eps$ such that 
\begin{enumerate}
\item[(i)] $F^\eps = \frac{1}{\eps} G^{\eps} + K^{\eps}$;
\item[(ii)] the Lie bracket of $G^\eps$ and $K^\eps$ vanishes, i.e. $[G^\eps,K^\eps]=0$; 
\item[(iii)] the vector field $G^\eps$ generates a  flow $\tau \mapsto \Phi^\eps_\tau$ which is $2\pi$-periodic, regardless of the specific trajectory, i.e.
$$
\forall (t,y) \in \R \times \R^n, \quad \Phi^\eps_{t+2 \pi}(y) = \Phi^\eps_{t}(y).
$$
\end{enumerate} 
\end{theorem}
This result brings us back to Section \ref{sect:split} and indeed allows to {\em split} equation (\ref{eq:transport}) into two equations of the form (\ref{eq:transporttau}-\ref{eq:transportt}); details will be given in Section \ref{sect:averkin}. We conclude this subsection with a few additional statements related to the conservation of geometric properties by stroboscopic averaging. 
\begin{theorem} \label{th:geo}
Suppose that the vector field $F^\eps$ can be split according to equation (\ref{eq:ode}) and that $G$ satisfies assumption $\bf (H)$. Then the two vector fields $G^\eps$ and $K^\eps$ of Theorem \ref{th:aver} have the following properties:
\begin{enumerate}
\item[(i)] if both $G$ and $K$ are divergence-free vector fields, then so are $G^\eps$ and $K^\eps$;
\item[(ii)] if both $G$ and $K$ are Hamiltonian vector fields, then so are $G^\eps$ and $K^\eps$.
\end{enumerate}
\end{theorem}
\begin{remark}
The properties of Theorem \ref{th:geo} are intimately linked to the choice of stroboscopic averaging (see \cite{chartier2,chartier1}), which is the only  averaging procedure preserving geometric properties of the initial vector field $F^\eps$.
\end{remark}
\subsection{Expansions in powers of $\eps$ of the vector fields $G^\eps$ and $K^\eps$} \label{sect:expansions}
Since we wish in particular to identify the asymptotic behaviour of (\ref{eq:transport}) in the limit when $\eps$ tends to zero as well as higher-order terms in $\eps$, it is essential to consider $\eps$-expansions of the various functions appearing in Theorem \ref{th:aver}. Since this was precisely the point of view adopted in \cite{chartier2,chartier1}, we shall again quote the following result\footnote{Note again that an alternative proof of this result may be found in \cite{murua1} and \cite{murua2}.}: 
\begin{theorem} \label{th:series}
Consider the Fourier series  of 
\begin{align} \label{eq:fourierK}
K_\tau (y) &= \left(\frac{\partial \Phi_{\tau}}{\partial y}(y)  \right)^{-1} \; \;  (K \circ \Phi_\tau) (y) 
= \sum_{k \in \Z} e^{i k \tau} \hat K_k(y).
\end{align}
The {\em averaged} vector field $K^\eps$ admits the following {\bf formal} $\eps$-expansion 
\begin{align} \label{eq:Kseries}
K^\eps &= \sum_{r=1}^{+\infty} \eps^{r-1} K^{[r]} 
       = \sum_{r=1}^{+\infty} \frac{\eps^{r-1}}{r} \sum_{(i_1, \ldots,i_r) \in \Z^r} \bar \beta_{i_1 \cdots i_r} \;  [\ldots[\hat K_{i_1},\hat K_{i_2}],\hat K_{i_3}],\ldots,\hat K_{i_r}]
\end{align}
where the coefficients $\bar \beta$ are universal (problem-independent). Similarly, the vector field $G^\eps$ admits the following formal $\eps$-expansion 
\begin{eqnarray} \label{eq:Gseries}
G^\eps &=& \eps (F^\eps - K^\eps).
\end{eqnarray}
\end{theorem}
\begin{remark}
The fact that geometric properties of $G^\eps$ and $K^\eps$ are inherited from $F^\eps$ may also be seen as a direct consequence of the form of previous expansions, which are linear combinations of embedded Lie-brackets of the $\hat K_k$'s. 
For instance, if both $G$ and $K$ are Hamiltonian, then $K_\tau$ is of the form 
$$
K_\tau(y) = J^{-1} \nabla_y H_\tau(y) \quad \mbox{ with } \quad  H_\tau(y) = \sum_{k \in \Z} e^{i \tau} \hat H_k(y)
$$
and  all Fourier coefficients $\hat K_k(y)= J^{-1} \nabla_y \hat H_k(y)$ are also Hamiltonian. Since 
$$
\forall (k,l) \in \Z^2, \quad [\hat K_k, \hat K_l] = J^{-1} \nabla_y \{\hat H_k, \hat H_l\}
$$ 
where $\{\cdot, \cdot\}$ denotes the Poisson bracket operation, it is then immediate to see that both $G^\eps$ and $K^\eps$ are Hamiltonian with Hamiltonians given by formulas (\ref{eq:Kseries}) and (\ref{eq:Gseries}) where Lie brackets are replaced by Poisson brackets and the $\hat K_k$'s by the $\hat H_k$'s. Similarly, if $\diver(G)=\diver(K)=0$, then $\diver(\hat K_k)=0$ for all $k \in \Z$ and a standard computation shows that 
$$
\forall (k,l) \in \Z^2, \quad \diver \left( [\hat K_k,\hat K_l] \right)=0
$$
so that again both $G^\eps$ and $K^\eps$ are divergence-free.
\end{remark}
In order to be able to derive the expansions of $G^\eps$ and $K^\eps$, it still remains to give the value of the coefficients $\bar \beta$ appearing in formula (\ref{eq:Kseries}). This is the purpose of next proposition.
\begin{proposition}
The coefficients $\bar \beta$ can be computed recursively from the following formulas, which hold for all values of $j \in \Z^{*}$, $r,s  \in \N^{*}$ and 
$(l_1, \ldots,l_s)  \in \Z^s$:
\begin{eqnarray*}
\begin{array}{llllll}
\bar \beta_0 &=& 1, & \bar \beta_j &=& 0, \\
\bar \beta_{0^{r+1}} &=& 0, &
\bar \beta_{0^{r}j} &=& \frac{i}{j} \left(\bar \beta_{0^{r-1}j} - \bar \beta_{0^r}\right), \\
\bar \beta_{j l_1 \cdots l_s} &=& \frac{i}{j} \left(\bar \beta_{l_1 \cdots l_s} - \bar \beta_{(j+l_1) l_2 \cdots l_s}\right), &
\bar \beta_{0^r j l_1 \cdots l_s} &=& \frac{i}{j} \left(\bar \beta_{0^{r-1} j l_1 \cdots l_s} - \bar \beta_{0^r (j+l_1) l_2 \cdots l_s}\right). 
\end{array}
\end{eqnarray*}
\end{proposition}
For the sake of illustration and later use, we now give the first terms of $K^\eps= K^{[1]} + \eps K^{[2]} + \eps^2 K^{[3]}  + {\cal O}(\eps^3)$, as stated in \cite{chartier4}:
\begin{align} \label{eq:firstK}
K^{[1]} &= \hat K_0, \nonumber \\
K^{[2]} &= \sum_{k >0} \frac{i}{k} \left( [\hat K_k, \hat K_{-k}] + [\hat K_0, \hat K_k-\hat K_{-k}]\right), \nonumber \\
K^{[3]} &= \sum_{k \neq 0} \frac{1}{k^2} \left( [[\hat K_k, \hat K_0],\hat K_0]+[[\hat K_{-k}, \hat K_k],\hat K_k]-\frac12 [[\hat K_{-2k}, \hat K_k],\hat K_k]+[ [\hat K_{0}, \hat K_k],\hat K_{-k}]\right) \nonumber \\
&- \sum_{0 \neq m \neq -l \neq 0} \frac{1}{l (m+l)} [[\hat K_{0}, \hat K_l],\hat K_{m}]
+ \sum_{k < -|l|} \frac{1}{l k} [[\hat K_{k}, \hat K_l],\hat K_{-l}]   \nonumber \\
&- \sum_{0 > k < m, m+k \neq 0} \frac{1}{k m } [[\hat K_{k}, \hat K_{-k}],\hat K_{m}] \nonumber \\
&- \sum_{0 \neq m \neq \pm l\neq 0, m >-m-l < l} \frac{1}{m (m+l)} [[\hat K_{-m-l},\hat K_{l}],\hat K_{m}].
\end{align}
\begin{remark}
The following expressions of the first three terms of the averaged equation have also been derived in various places and do not use Fourier coefficients:
\begin{align*}
K^{[1]}(y) &= \frac{1}{2 \pi} \int_0^{2 \pi} K_\tau(y) d\tau, \quad
K^{[2]}(y) = \frac{-1}{4 \pi} \int_0^{2 \pi}  \hskip -1ex \int_0^\tau [K_s(y),K_\tau(y)] ds d\tau,  \\
K^{[3]}(y) &= \frac{1}{8 \pi} \int_0^{2 \pi}  \hskip -1ex \int_0^\tau \hskip -1ex \int_0^s [[K_r(y),K_s(y)],K_\tau(y)] dr ds d\tau \\
& \qquad \qquad + \frac{1}{24 \pi} \int_0^{2 \pi}  \hskip -1ex \int_0^\tau \hskip -1ex \int_0^\tau [K_r(y),[K_s(y),K_\tau(y)]] dr ds d\tau.
\end{align*}
 Further terms can be formally obtained  by using a non-linear Magnus expansion \cite{blanes}.  Each of these is a linear combination of  iterated integrals of iterated brackets of $K_\tau$.
\end{remark}
As an illustration, we derive below the expressions of $G^\eps$ and $K^\eps$ for a simple example. We thus consider the following vector field 
\begin{align} \label{eq:Fex}
F^\eps(y) = 
\left(
\begin{array}{c}
v \\
\frac{1}{\eps} J v + E
\end{array}
\right)
\end{align}
where $x=(x_1,x_2)^T \in \R^2$, $v=(v_1,v_2)^T \in \R^2$, 
$y=(x_1,x_2,v_1,v_2)^T \in \R^4$, $E\in \R^2$ and 
$$
J=\left(
\begin{array}{cc}
0 & 1 \\
-1 & 0
\end{array}
\right).
$$ 
The function $F^\eps$ may be decomposed into the sum $\frac{1}{\eps} G + K$ with 
$$
G(y) = 
\left(
\begin{array}{c}
0 \\
J v 
\end{array}
\right)
\quad \mbox{ and } \quad 
K(y) = 
\left(
\begin{array}{c}
v \\
E 
\end{array}
\right),
$$
and the flow $\Phi_\tau$ associated with $G$  simply reads 
\begin{eqnarray*}
\Phi_\tau(y) = \left(
\begin{array}{c}
x\\
e^{\tau J} v
\end{array}
\right).
\end{eqnarray*}
Substituting $\Phi_\tau$ into $K$ then leads to 
\begin{eqnarray*}
K_\tau(y) = \left(
\begin{array}{c}
e^{\tau J} v\\
e^{-\tau J} E
\end{array}
\right) = e^{i \tau} \hat K_1(y) + e^{-i \tau} \hat K_{-1}(y) 
\end{eqnarray*}
with 
\begin{eqnarray*}
\hat K_1(y) = \frac12  \left(
\begin{array}{c}
v-i J v\\
E+iJ E
\end{array}
\right)
\quad \mbox{ and } \quad 
\hat K_{-1}(y) =
\frac12  \left(
\begin{array}{c}
v+ iJ v\\
E- iJ E
\end{array}
\right),
\end{eqnarray*}
where we have used the relation $e^{\theta J} = (\cos \theta) I + (\sin \theta) J$ and have written $\cos \theta=\frac12( e^{i \theta}+e^{-i\theta})$ and $\sin \theta=\frac{1}{2 i}( e^{i \theta}-e^{-i\theta})$. 
Formula (\ref{eq:Kseries}) then gives
\begin{eqnarray*}
K^{[1]} &=& \hat K_0 = 0, \\ 
K^{[2]} &=& i [\hat K_1,\hat K_{-1}] = -2 \Im \left((\partial_y \hat K_1) \hat K_{-1} \right) = \left(
\begin{array}{c}
J E\\
0
\end{array}
\right)
\end{eqnarray*}
and all other $K^{[r]}$ for $r \geq 3$ vanish, as can be checked by easy calculations.
\section{Averaging of transport equations with constant frequency} \label{sect:averkin}
Compiling the arguments of the two previous sections, it is now straightforward to obtain the following corollary, which establishes in particular the  existence of a  formal averaged transport equation for problems of the form (\ref{eq:transport}-\ref{eq:splitF}). 
\begin{corollary} \label{th:main}
Let $F^\eps=\frac{1}{\eps} G^\eps + K^\eps$ be the normal form splitting of a highly-oscillating vector field $F^\eps=\frac{1}{\eps} G + K$ satisfying $\bf (H)$. The solution of the transport equation 
$$
\partial_t f(t,y) + F^\eps(y) \cdot \nabla_y f(t,y)=0
$$
may be obtained as the diagonal value (i.e. for the value $\tau=t/\eps$) of the two-scale function $\tilde f(t,\tau,y)$, $2 \pi$-periodic in $\tau$, defined as the unique solution of the following system of two equations
\begin{eqnarray*} 
\left\{
\begin{array}{c}
   \forall (t,\tau,y), \quad \partial_\tau \tilde f(t,\tau,y) + G^\eps(y) \cdot \nabla_y \tilde f (t,\tau,y) =  0, \quad { (i)}\\
 \forall (t,\tau,y), \quad \partial_t \tilde f(t,\tau,y) + K^\eps(y) \cdot \nabla_y \tilde f (t,\tau,y) = 0 \quad { (ii)} 
\end{array}
\right.
\end{eqnarray*}
with initial condition $\tilde f(0,0,\cdot)=f_0$. 
Moreover, the $\eps$-expansions of $G^\eps$ and $K^\eps$ are given by formulas (\ref{eq:Kseries}-\ref{eq:Gseries}) of Theorem \ref{th:series}.     If in addition $G$ and $K$ are both divergence-free, then so are $G^\eps$ and $K^\eps$, and similarly, if $G$ and $K$ are both Hamiltonian, then so are $G^\eps$ and $K^\eps$, with Hamiltonians that can be obtained again from formulas (\ref{eq:Kseries}-\ref{eq:Gseries}) by replacing Lie brackets by Poisson brackets. 
\end{corollary}
\begin{Proof} The result follows 
immediately from Proposition \ref{lem:splittrans} with $F_1 = \frac{1}{\eps} G^\eps$ and $F_2 = K^\eps$ and from Theorem \ref{th:series}. 
\end{Proof}
\begin{remark} Equation $(ii)$ is usually referred to as  the {\em averaged transport equation}.
\end{remark}
As a straightforward illustration of this corollary, we consider the simplified case of a set of particles evolving in a constant electric field (independent of time and phase-space variables) and submitted to a constant magnetic field. The corresponding equation 
\begin{eqnarray} \label{eq:transportex}
\partial_t f + v \cdot \nabla_x f + \left( \frac{1}{\eps} J v + E \right) \cdot \nabla_v f = 0,
\end{eqnarray}
-where $f$ depends on time $t \in \R$, position $x \in \R^2$ and velocity $v \in \R^2$- is obviously of the form (\ref{eq:transport}) with $y=(x_1,x_2,v_1,v_2)^T \in \R^4$ and $F^\eps$ given by (\ref{eq:Fex}). 
On the one hand, given the extreme simplicity of the vector field $F^\eps$, the solution $f(t,x,v)$ of (\ref{eq:transportex}) can be directly written  as 
\begin{eqnarray} \label{eq:exact}
f_0\left(x-\eps J(e^{\tau J}-I) v-\eps^2 (e^{\tau J}E-E)+\eps t J E,e^{\tau J}v - \eps J(e^{\tau J}-I)E\right)
\end{eqnarray}
for $\tau=t/\eps$. On the other hand, using the computations at the end of previous section, equations $(i)$ and $(ii)$ of Corollary \ref{th:main} for $\tilde f(t,\tau,x,v)$  take the following form
\begin{eqnarray*}
 (i) \quad \partial_\tau \tilde f + (\eps v - \eps^2 JE) \cdot \nabla_x \tilde f + (Jv + \eps E) \cdot \nabla_v \tilde f = 0, \qquad
 (ii) \quad \partial_t \tilde f + \eps J E \cdot \nabla_x \tilde f= 0.
\end{eqnarray*}
By direct differentiation w.r.t. $\tau$ and then $t$, it can be checked that the function given in formula (\ref{eq:exact})  satisfies both equations $(i)$ and $(ii)$.  

\section{High-oscillations with varying frequency} \label{sect:varyfreq}
In this section, we again consider the transport equation  
\begin{eqnarray} \label{eq:checktransport}
\partial_t f(t,y) + F^\eps(y) \cdot \nabla_y f(t,y) = 0
\end{eqnarray}
where the vector field  $F^\eps$ is now of the form
\begin{eqnarray} \label{eq:checkF}
F^\eps(y) = \frac{1}{\eps}  \omega(y) G(y) +  K(y)
\end{eqnarray}
with $G$ still generating a $2 \pi$-periodic flow 
$\Phi_\tau$, independently of the initial condition. In this form, Theorem \ref{th:main} does not directly apply, owing to the non-existence of a common frequency for all trajectories (if $\omega$ varies). In order to rewrite (\ref{eq:checktransport}) in a more amenable form, we thus divide it by $\omega$
\begin{eqnarray} \label{eq:Lie}
\frac{1}{\omega(y)} \partial_t f(t,y) + \frac{1}{\omega(y)} F^\eps(y) \cdot \nabla_y f(t,y) = 0.
\end{eqnarray}
Upon denoting $Y=(t,y)$, previous equation may then be rewritten as ${\mathcal L}_{\check F^\eps} (f) =0$, where 
\begin{eqnarray} \label{eq:Lie2}
 {\mathcal L}_{\check F^\eps} (f)=  (\pa_Y f) \; \check F^\eps
\end{eqnarray}
 is the Lie derivative of $f$ in the direction of the {\bf augmented} vector field 
\begin{eqnarray} \label{eq:augF}
\check F^\eps(Y)&=&
\left(
\begin{array}{c}
\frac{1}{\omega(y)} \\
\frac{1}{\omega(y)} F^\eps(y)
\end{array}
\right)=
\frac{1}{\eps}
\left(
\begin{array}{c}
0 \nonumber \\
G(y)
\end{array}
\right)+
\left(
\begin{array}{c}
\frac{1}{\omega(y)} \\
\frac{1}{\omega(y)} K(y) 
\end{array} 
\right)\\
&:=& \frac{1}{\eps} \check G(Y) + \check K(Y).
\end{eqnarray}
In particular, note that $\check G$ still generates a $2 \pi$-periodic flow. 

\subsection{Immersion as the stationary solution of an extended equation}
Our first idea  is to interpret the  function $f(t,y)=f(Y)$  as the (stationary) solution to  the following augmented transport equation  on $g(s,Y)$:
$$
\partial_s g(s,Y) + \check F^\eps(Y) \cdot \nabla_Y g(s,Y) = 0, \quad g(0,Y)=f(Y)=f(t,y).
$$
This means that 
$$
g(s,\cdot) = \exp\left(-s {\mathcal L}_{\check F^\eps} \right) f = f \quad \mbox{for all}  \quad s\geq 0.
$$
 Since $\check G$  generates a $2 \pi$-periodic flow, the averaging Theorem \ref{th:aver} ensures that
 $$\check F^\eps = \frac{1}{\eps} \check G^\eps + \check K^\eps,
 $$
 where $\check G^\eps$ still generates a $2 \pi$-periodic flow and  $[\check G^\eps, \check K^\eps]=0$.
Proceeding as in Section \ref{sect:split}, we then get two equations for 
$$
\tilde  g(s,\tau,\cdot) = \exp\left(-\tau {\mathcal L}_{\check G^\eps} \right)\exp\left(-s {\mathcal L}_{\check K^\eps} \right)  f
$$
of the form
\begin{eqnarray} 
(i) & \partial_s \tilde g(s,\tau,Y) + \check K^\eps(Y) \cdot \nabla_Y  \tilde g(s,\tau,Y) = 0 \label{two-scale-moy1} \\
(ii) & \partial_\tau   \tilde g(s,\tau,Y) + \check G^\eps(Y) \cdot \nabla_Y \tilde g(s,\tau,Y) = 0 \label{two-scale-moy2}
\end{eqnarray}
which  can be solved one after another in any order, since $[\check G^\eps, \check K^\eps]=0$. Note the usual relation $\tilde g(s,s/\eps,Y) = g(s,Y)=f(Y)$. However, there is here no known  initial condition at $s=\tau=0$, since $\tilde g (0,0,Y)=g(0,Y)=f(t,y)$ is precisely the unknown of the original problem.

  \subsection{Eliminating the extra-variable $s$}
  
  Our objective in this subsection is to transform the two equations (\ref{two-scale-moy1}-\ref{two-scale-moy2})  into new equations which do not involve the variable $s$ and are provided with a proper initial condition, namely $f_0(y)$. We will then show how to recover the original solution $f(t,y)$ using only these new equations. To this aim, we will introduce a phase-function $(t,\tau,y) \mapsto S(t,\tau,y)$ in the spirit of \cite{crouseilles}, which will be defined later on as the solution of a transport equation,  and a profile-function $(t,\tau,y) \mapsto h(t,\tau,y)$ defined by 
\begin{eqnarray}
\label{link-hg}
h(t,\tau, y)= \tilde g ( S(t,\tau,y),\tau,t,y),
\end{eqnarray}
that will also be shown to satisfy a companion transport equation. Our starting point is the following set of relations
\begin{eqnarray*}
\pa_t h &=& \left(\pa_s \tilde g( S,\tau,t,y) \right) \pa_t S + \pa_t \tilde g( S,\tau,t,y), \\
\pa_\tau h &=& \left(\pa_s \tilde g( S,\tau,t,y) \right) \pa_\tau S + \pa_\tau\tilde g( S,\tau,t,y),\\
\pa_y h&=& \left(\pa_s \tilde g( S,\tau,t,y) \right) \pa_y S + \pa_y \tilde g( S,\tau,t,y), 
\end{eqnarray*}
where we have omitted the obvious arguments of functions $h$ and $S$ and which may be straightforwardly obtained. Together with equations \fref{two-scale-moy1} and \fref{two-scale-moy2}, they lead immediately to
\begin{eqnarray*}
&\check K_1^\eps(y) \pa_t h(t,\tau, y) + \check K_2^\eps(y) \cdot \nabla_y h(t,\tau, y) \\
&= \left(\pa_s \tilde g( S,\tau,t,y) \right) \left( \check K_1^\eps(y) \pa_t S(t,\tau, y) + \check K_2^\eps(y) \pa_y S(t,\tau, y) -1\right),
\end{eqnarray*}
  and
\begin{eqnarray*}
&\pa_\tau h(t,\tau, y)  + \check G_1^\eps(y) \pa_t h(t,\tau, y) + \check G_2^\eps(y) \cdot \nabla_y h(t,\tau, y) \\
&=  \left(\pa_s \tilde g( S(t,\tau,y),\tau,t,y) \right) \left( \check G_1^\eps(y) \pa_t S(t,\tau, y) + \check G_2^\eps(y) \pa_y S(t,\tau, y) +\pa_\tau S(t,\tau, y)\right),
\end{eqnarray*}
where the index $1$ in $\check K_1^\eps$ and $\check G_1^\eps$ refers to the first components of $\check K^\eps$ and $\check G^\eps$, while the index $2$ in $\check K_2^\eps$ and $\check G_2^\eps$ refers to {\em all remaining components} of $\check K^\eps$ and $\check G^\eps$. 
Now, in order to eliminate the variable $s$ from the previous two equations, one has to choose $S$  such that
 \be
 \label{eq:S}
 \begin{array}{l}
 \check K_1^\eps(y) \pa_t S(t,\tau, y) +\check K_2^\eps(y) \cdot \nabla_y S(t,\tau, y) = 1, \\
 \pa_\tau S(t,\tau, y)+ \check G_1^\eps(y) \pa_t S(t,\tau, y) + \check G_2^\eps(y) \cdot \nabla_y S(t,\tau, y)=0,
 \end{array}
 \ee
 and then
 \be
 \label{eq:h}
 \begin{array}{l}
 \check K_1^\eps(y) \pa_t h(t,\tau, y) + \check K_2^\eps(y) \cdot  \nabla_y h(t,\tau, y) =0, \\
 \pa_\tau h(t,\tau, y) + \check G_1^\eps(y) \pa_t h(t,\tau, y) + \check G_2^\eps(y) \cdot \nabla_y h(t,\tau, y)=0,
 \end{array}
 \ee
 with initial conditions  
\begin{eqnarray} \label{eq:initcond}
S(0,0,y)=0, \qquad h(0,0,y)= f_0(y).
 \end{eqnarray}
 From these functions $S$ and $h$, one can recover the distribution function $f(t,y)$ as follows: for any given $(t,y)$, define $\tau(t,y)$ as a solution of
  $$\tau(t,y)= \frac{S(t,\tau(t,y),y)}{\eps}.$$
  Then $f$ can be obtained from the relation
  \begin{align*}
  h(t, \tau(t,y),y)&= \tilde g \left( S(t,\tau(t,y),y),  \frac{S(t,\tau(t,y),y)}{\eps} , t,y\right),\\
  &= g (S(t,\tau(t,y),y),t,y)\\
  &=f(t,y).
  \end{align*}
%
%
\begin{lemma} \label{lem:commu}
Assume that $y \mapsto \check K_1^\eps(y)$ does not vanish and consider the two vector fields
\begin{align*}
\check A^\eps := \frac{1}{\check K_1^\eps} \, \check K_2^\eps \quad \mbox{ and } \quad \check B^\eps :=\check G_2^\eps - \frac{\check G_1^\eps}{\check K_1^\eps} \, \check K_2^\eps,
\end{align*}
together with the two scalar functions
\begin{align*}
\check \alpha^\eps := \frac{1}{\check K_1^\eps}  \quad \mbox{ and } \quad \check \beta^\eps :=- \frac{\check G_1^\eps}{\check K_1^\eps}.
\end{align*}
Then the following two relations hold true
\begin{align}
{\cal L}_{\check A} \check \beta = {\cal L}_{\check B} \check \alpha  \quad \mbox{ and } \quad {\cal L}_{\check A} {\cal L}_{\check B} = {\cal L}_{\check B} {\cal L}_{\check A}. 
\end{align}
\end{lemma}  
\begin{Proof}
Owing to Theorem \ref{th:series}, the two vector fields $\check K^\eps$ and $\check G^\eps$ 
have a vanishing Lie bracket (with respect to the $Y=(t, y)$ variable). This implies that 
\begin{align} \label{eq:comu}
 \partial_y \check K_1^\eps(y) \; \check G_2^\eps(y) -  \partial_y \check G_1^\eps(y) \; \check K_2^\eps(y)=0 \; \mbox{ and } \; 
\partial_y \check K_2^\eps(y) \; \check G_2^\eps(y) -  \partial_y \check G_2^\eps(y) \; \check K_2^\eps(y)=0.
\end{align}
By definition of $\check \alpha^{\eps}$ and $\check \beta^{\eps}$, the first relation may be rewritten as 
\begin{align} \label{eq:part2}
\nabla_y \check \beta^{\eps} \cdot \check A^\eps - \nabla_y \check \alpha^{\eps} \cdot \check B^\eps = 0
\end{align}
which proves the first statement of the lemma. Now, given a smooth vector field $L: \R^{n} \rightarrow \R^{n}$, a scalar function $a:\R^{n} \rightarrow \R$ and $\delta y$ a vector of $\R^{n}$, the relation 
$$
\left( \partial_y ( a L ) \right) \, \delta y = \left(\nabla_y a \cdot \delta y \right) L + a (\partial_y L) \, \delta y
$$
holds true and may be used to compute the Lie bracket of $\check A^{\eps}$ and $\check B^{\eps}$ as follows
\begin{align*}
[\check A^{\eps},\check B^{\eps}] =&  \partial_y ( \check \alpha^{\eps} \check K_2^{\eps} )  \, \left( \check G_2^\eps +\check \beta^\eps \check K_2^\eps \right) - \partial_y \check G_2^\eps \, ( \check \alpha^{\eps} \check K_2^{\eps} ) - \partial_y ( \check \beta^{\eps} \check K_2^{\eps} ) \, ( \check \alpha^{\eps} \check K_2^{\eps} ) \\
=&\left(\nabla_y \check \alpha^\eps \cdot \left( \check G_2^\eps +\check \beta^\eps  \check K_2^\eps \right) \right) \check K_2^{\eps} + \check \alpha^{\eps} (\partial_y \check K_2^\eps) \, \left( \check G_2^\eps +\check \beta^\eps \check K_2^\eps \right) \\
& - \check \alpha^{\eps} (\partial_y \check G_2^\eps) \,  \check K_2^{\eps} - \check \alpha^{\eps} \left(\nabla_y \check \beta^{\eps} \cdot  \check K_2^{\eps}  \right) \check K_2^{\eps} -\check \alpha^{\eps} \check \beta^{\eps} (\partial_y \check K_2^{\eps}) \,  \check K_2^{\eps}.
\end{align*}
Using the second half of \fref{eq:comu}, the equality above simplifies to 
\begin{align*}
[\check A^{\eps},\check B^{\eps}] =&\left(\nabla_y \check \alpha^\eps \cdot  \check G_2^\eps + \check \beta^\eps  \nabla_y \check \alpha^\eps \cdot  \check K_2^\eps - \check \alpha^{\eps} \nabla_y \check \beta^{\eps} \cdot  \check K_2^{\eps}  \right) \check K_2^{\eps} \\
=&\left( \nabla_y \check \alpha^\eps \cdot  \check B^\eps - \nabla_y \check \beta^{\eps} \cdot  \check A^{\eps} \right) \check K_2^{\eps}
\end{align*}
where the scalar term  in factor of $\check K^{\eps}_2$ now vanishes owing to \fref{eq:part2}. This implies the second statement of the lemma and completes its proof.
\end{Proof}

  %
  %
\begin{theorem}
\label{th-S}
Consider the functions  $S(t,\tau,y)$ and $h(t,\tau,y)$ satisfying the following two separate  systems of equations
 \begin{align}
 & \check K_1^\eps(y) \pa_t S(t,\tau, y) + \check K_2^\eps(y) \cdot \nabla_y S(t,\tau, y) = 1, \label{eq:St} \\
 &\check K_1^\eps(y) \pa_\tau S(t,\tau, y)+ \Big(\check K_1^\eps(y) \check G_2^\eps(y) - \check G_1^\eps(y) \check K_2^\eps(y) \Big) \cdot \nabla_y S(t,\tau, y)=-\check G_1^\eps(y), \label{eq:Stau}\\
 &S(0,0,y)=0, \label{eq:Sinit}
 \end{align}
 and
 \begin{align}
 & \check K_1^\eps(y) \pa_t h(t,\tau, y) + \check K_2^\eps(y) \cdot  \nabla_y h(t,\tau, y) =0, \label{eq:ht} \\
 & \check K_1^\eps(y) \pa_\tau h(t,\tau, y)+ \Big(\check K_1^\eps(y) \check G_2^\eps(y) - \check G_1^\eps(y) \check K_2^\eps(y) \Big) \cdot \nabla_y h(t,\tau, y)=0, \label{eq:htau}\\
 & h(0,0,y)= f_0(y). \label{eq:hinit}
 \end{align}
If the function $y \mapsto \check K_1^\eps(y)$ does not vanish, then the following statements hold:
 \begin{enumerate}
\item[(i)]  system  (\ref{eq:ht}-\ref{eq:htau}-\ref{eq:hinit}) has a unique solution $h$, periodic w.r.t. $\tau$;
\item[(ii)] system (\ref{eq:St}-\ref{eq:Stau}-\ref{eq:Sinit}) has a unique solution $S$, periodic w.r.t. $\tau$; 
\item[(iii)] the  {\bf \em formal} expansion of the solution $f(t,y)$ of problem (\ref{eq:checktransport}-\ref{eq:checkF}) satisfies
$$
f(t,y) = h(t,\tau(t,y),y),  
$$
where the function $(t,y) \mapsto \tau(t,y) \in \R$ is implicitly defined (locally) by the relation 
$$
\eps \tau(t,y)= S(t,\tau(t,y),y).
$$
\end{enumerate}
\end{theorem}
\begin{Proof}
A straightforward computation shows that the four equations \fref{eq:St}, \fref{eq:Stau}, \fref{eq:ht}, \fref{eq:htau} are equivalent to the four equations in (\ref{eq:S}) and (\ref{eq:h}). Hence, if the separate systems (\ref{eq:ht}-\ref{eq:htau}) and (\ref{eq:St}-\ref{eq:Stau}) have unique solutions, they are clearly periodic w.r.t. $\tau$. Now, proving the first statement requires to show that equations \fref{eq:ht} and \fref{eq:htau} can be solved in any order, i.e. that ${\cal L}_{\check A^\eps}$ and ${\cal L}_{\check B^\eps}$ commute, which is ensured by previous lemma. If $h$ is the solution of (\ref{eq:ht}-\ref{eq:htau}-\ref{eq:hinit}), then $h(\cdot,0,\cdot)$ is in particular the solution of the Cauchy problem (\ref{eq:ht}-\ref{eq:hinit}) and thus reads
$$ 
h(t,0,\cdot)=  \exp(-t\mathcal L _{\check A^\eps})f_0.
$$
Equation \fref{eq:htau}, which is a transport equation in variables $(\tau,y)$ with fixed parameter $t$, can then be uniquely solved. Given the initial data $h(t,0,\cdot)=  \exp(-t\mathcal L _{\check A^\eps})f_0$, this yields
\begin{align} \label{eq:httau}
h(t,\tau,\cdot)=  \exp(-\tau \mathcal L _{\check B^\eps}) \exp(-t\mathcal L _{\check A^\eps})f_0.
\end{align}
Hence, if a solution of  (\ref{eq:ht}-\ref{eq:htau}-\ref{eq:hinit}) exists, it is necessarily of this form and thus unique. Conversely, one has, according to previous lemma 
$$
 \exp(-t\mathcal L _{\check A^\eps}) \exp(-\tau \mathcal L _{\check B^\eps}) f_0 =\exp(-\tau \mathcal L _{\check B^\eps}) \exp(-t\mathcal L _{\check A^\eps})f_0
$$
and by differentiating the left-hand side w.r.t. $t$ and the right-hand side $\tau$, it may be checked that $h$ given in \fref{eq:httau}  is indeed solution -thus the unique solution- of  system (\ref{eq:ht}-\ref{eq:htau}-\ref{eq:hinit}). This proves (i).

Proceeding similarly for system (\ref{eq:St}-\ref{eq:Stau}-\ref{eq:Sinit}), we first solve (\ref{eq:St}-\ref{eq:Sinit}) for fixed $\tau =0$. This yields
$$
S(t,0,\cdot) = \exp\big(-t {\cal L}_{\check A^\eps} \big) \, S(0,0,\cdot) +  \int_0^t \exp\big((s-t) {\cal L}_{\check A^\eps}\big) \, ds \; \alpha = t \varphi\big(-t {\cal L}_{\check A^\eps} \big) \, \alpha,
$$
where $\varphi(z) = \frac{e^z-1}{z}$ is holomorphic on $\C$.
The function $S$ so-obtained then serves as initial  condition for the evolution in $\tau$ through  equation \fref{eq:Stau}. This then leads to 
\begin{align*}
S(t,0,\cdot) &= \exp\big(-\tau {\cal L}_{\check B^\eps} \big) t \varphi\big(-t {\cal L}_{\check A^\eps} \big) \, \alpha + \tau \varphi\big(-\tau {\cal L}_{\check B^\eps} \big) \, \beta \\
&=\tau \varphi\big(-\tau {\cal L}_{\check B^\eps} \big) \, \beta + t \varphi\big(-t {\cal L}_{\check A^\eps} \big) \alpha -t \, \tau \,  \varphi\big(-t {\cal L}_{\check A^\eps} \big) \varphi\big(-\tau {\cal L}_{\check B^\eps} \big) {\cal L}_{{\check B^\eps}} \alpha
\end{align*}
where we have used the commutation of ${\cal L}_{\check A^\eps}$ and ${\cal L}_{\check B^\eps}$. Solving the equations in reverse order would have led to the symmetric variant
$$
S(t,0,\cdot) = \tau \varphi\big(-\tau {\cal L}_{\check B^\eps} \big) \, \beta + t \varphi\big(-t {\cal L}_{\check A^\eps} \big) \alpha -t \, \tau \,  \varphi\big(-t {\cal L}_{\check A^\eps} \big) \varphi\big(-\tau {\cal L}_{\check B^\eps} \big) {\cal L}_{\check A^\eps} \beta,
$$
which, owing to Lemma \ref{lem:commu} (${\cal L}_{\check A^\eps} \beta = {\cal L}_{\check B^\eps} \alpha$), coincides with the first one. This proves {\em (ii)}. 

It remains to prove {\em (iii)}. From \fref{link-hg} and the definition of $h$ and 
$\tilde g$, we have 
$$ 
  \forall \, (t,\tau ,y), \quad h(t, \tau ,y)= \tilde g \left( S(t,\tau,y), \tau , t,y\right) \quad \mbox{ and } \quad 
\tilde g \left( S(t,\tau,y), \frac{S(t,\tau,y)}{\eps} , t,y\right) = f(t,y), 
$$ 
so that the value of $f(t,y)$  can be recovered from $h$ and $S$  through the formula 
$$
\forall \, (t,\tau ,y), \quad f(t,y) = h(t,\tau(t,y),y),
$$
provided that  $\tau(t,y)$ satisfies
$$
\eps \tau(t,y) = S(t,\tau(t,y),y).
$$
Given the periodicity of $S$ w.r.t. $\tau$, this equation always has a solution $\tau(t,y)$.
\end{Proof}

 \begin{remark}{\bf (truncated averaged models)} 
 \label{remark-th-S}
 If one keeps, in the expansions of the averaged fields ${\check A^\eps}$ and ${\check B^\eps}$ (defined in Lemma \ref{lem:commu}), only the terms of order less than (or equal to) $n$ in $\eps$, then the question arises whether the corresponding  truncated averaged models\footnote{i.e. the models obtained by removing all the terms of size $\eps^p$ for $p\geq n+1$.}  have a solution in the exact sense, and whether this solution is periodic w.r.t. $\tau$.  Generally speaking, the transport operators associated with  the truncated fields 
 $\check A_n$ and $\check B_n$ 
 (i.e. $\check A^\eps=\check A_n+O(\eps^{n+1})$,  
 $\check B^\eps=\check B_n+O(\eps^{n+1})$) do not commute exactly. More precisely, we only have $[\check A_n, \check B_n]= O(\eps^{n+1})$. However, one can define an approximate solution by first solving
$$ 
\pa_\tau h + \check B_n(y)\cdot \nabla_y h=0, \quad h(0,0,y)=f_0(y)
$$ 
for fixed $t=0$ (in this way we obtain a solution $h_n(0,\tau,y)$ defined for all $\tau$), and then solving 
$$ 
\pa_t h + \check  A_n(y)\cdot \nabla_y h=0, \quad h(0,\tau,y)=h_n(0,\tau,y)
$$
in order to get a solution $h_n^1(t,\tau,y)$ defined for all $\tau$  and $t$.  At this stage, it is worth emphasizing that the function $h^1_n$ does not  satisfy exactly the first equation for all $t$ (only for $t=0$), since $[\check A_n, \check B_n]\ne 0.$  Nevertheless, it  does satisfy it up to terms
 of size $\eps^{n+1}$.  In particular, if one solves the two equations in reverse order, the function $h^2_n$ obtained does not coincide with $h_n^1$ exactly, but only up to terms of size $\eps^{n}$ and we have $h_n^1-h_n^2 = O(\eps^{n+1})$. 
 In this sense, the result in previous theorem is at this stage only {\em formal}. It will be the subject of a forthcoming paper \cite{chartierfuture} to prove error estimates for the defects in (\ref{eq:S}) and (\ref{eq:h}).
 \end{remark}
%

\subsection{An illustrative elementary example}
Our aim here is to illustrate the result of Section \ref{sect:varyfreq} on an elementary example for which exact solutions can be easily obtained. Consider the following transport equation 
\be
\label{ill-example}
\pa_t f + \left(\frac{1}{\eps} \omega(y)Jy + y \right)\cdot \nabla_y f = 0, \qquad f(0,y)= f_0(y),
\ee
where $y \in \R^2$, and where
$$
 \omega (y) = 1+ |y|^2 = 1+ y_1^2 + y_2^2 \quad \mbox{ and } \quad J=\left( \begin{array}{rr} 0 & 1\\ -1 &0\end{array}\right) .$$
 This equation can be solved as follows: let $\varphi_t^\eps(y)$ be the flow of the characteristics equation
 $$
 \dot y  =  \frac{1}{\eps} \omega(y)Jy + y.
 $$
 By taking its inner product by $y$, we have immediately $|\varphi_t^\eps(y)| = \exp(t) |y|$, so that 
 $$
 \varphi_t^\eps(y)= \exp(t) \exp\left( \frac{1}{\eps} \left(t +    (e^{2t} -1) \frac{|y|^2}{2}\right)J \right)y.
 $$
 As a consequence, the explicit solution of \fref{ill-example} reads
 \be
 \label{sol-ill-example}
 f(t,y) = f_0 \left(\exp\left( -t  -\frac{1}{\eps} \left(t + (1- e^{-2t}) \frac{|y|^2}{2}\right)J \right)y\right) .
 \ee
 Now, we observe that the two fields $\omega(y) Jy$ and $K(y)= y$ do not commute, and in order to transform the problem into a highly-oscillatory
 problem  with $y$-independent frequency, one has to divide equation \fref{ill-example} by $\omega$ and immerse the 
 equation on $f$ into an augmented one for the unknown $g(s,t,y)$
 
 \be
 \label{imm-ill-problem}
 \pa_s g+ \frac{1}{\omega(y)} \pa_t g + \left(\frac{1}{\eps} Jy + \frac{y}{\omega(y)} \right)\cdot \nabla_y g = 0, \qquad g(0,t,y)= f(t,y).
 \ee
 Unlike the fields $\omega(y) Jy$ and $K$, we now observe that the two augmented fields  $\check  G(y) = (0,Jy)^T$ and $\check K(y) = \left(\frac{1}{\omega(y)}, \frac{y}{\omega(y)}\right)^T$ do  commute.  This means that equation \fref{imm-ill-problem} is already written in a normal form and therefore the averaged fields in this case are simply
 $$ \check G^\eps= (0,Jy)^T, \quad \check K^\eps=(\check K_1^\eps,\check K_2^\eps)^T, \quad \mbox{with} \quad \check K_1^\eps=\frac{1}{\omega(y)}, \quad \check K_2^\eps= 
  \frac{y}{\omega(y)}.$$
  We now apply  Theorem \ref{th-S} in this particular case. The solution $h=h(t,0,y)$ to
  $$\pa_t h + y\cdot \nabla_y h =0, \qquad  h(0,0,y)= f_0(y),$$ is $h(t,0,y)= f_0(e^{-t} y)$. As a consequence, the solution $h=h(t,\tau,y)$ to
$$\pa_\tau h + Jy \cdot \nabla_y h=0, \qquad  h(t,0,y)= f_0(e^{-t} y),
$$
is 
\be
h(t,\tau,y) = f_0(e^{-t} e^{-\tau J}  y).
\ee
  The solution $S=S(t,0,y)$ to 
  $$\pa_t S + y\cdot \nabla_y S=\omega(y), \qquad  S(0,0,y)= 0,$$ is simply $S(t,0,y)=  t + (1- e^{-2t}) \frac{|y|^2}{2}$, 
  so that the solution $S=S(t,\tau,y)$ to 
  $$\pa_\tau S + Jy \cdot \nabla_y S=0, \qquad  S(t,0,y)= t + (1- e^{-2t}) \frac{|y|^2}{2},$$
  is constant w.r.t. $\tau$, given that $|e^{\tau J} y|^2 = |y|^2$, i.e.
\be
S(t,\tau,y)=  t + (1- e^{-2t}) \frac{|y|^2}{2}.
\ee
Theorem \ref{th-S} asserts that  $f(t,y) = h(t,\tau(t,y),y)$ 
where $\tau(t,y)$ is given by $\eps \tau(t,y) = t+(1- e^{-2t}) \frac{|y|^2}{2}$, an assertion which can be easily checked on our explicit example.

\section{Application to Vlasov equations with a strong magnetic field} \label{sect:vlasov}
In this section, we consider the case of particles submitted to a strong magnetic field and evolving in an electric field $E(x)$ depending on the position $x$ only. We recall hereinafter the corresponding equation (\ref{eq:vlasov}) on the distribution function $ f=f(t,x,v)$,  $t\geq 0$, $x\in \R^3$, $v\in \R^3$:
\begin{equation} \label{eq:vlasovJ}
\partial_t f + v \cdot \nabla_x f+ \left(E(x) +v\times  \frac{B(x)}{\eps} \right) \cdot \nabla_v f = 0, \qquad f(0,x,v)=f_0(x,v),
\end{equation}
which is closely related to the illustrative example of Section \ref{sect:averkin}, though with the additional difficulty that $E$ and $B$ may vary. We further assume here that $E$ derives from a potential $U$, i.e. that $E(x) = -\nabla_x U(x)$. 
\subsection{Constant magnetic field} \label{sect:constant}
Over a first phase, we assume that the magnetic field is constant. This means that, up to constant rotation, we have $B(x) =(0,0, b(x))^T$ and $\ b(x) \equiv b$. Upon rescaling the time $t \rightarrow t/b$ in $f$, i.e. considering the equation for $f(t/b,x,v)$ instead of $f(t,x,v)$ we may even assume that $b=1$.  We further assume in this first phase that the potential $U$ depends only on the orthogonal direction (to $B$) of $x$, that is on the first two  components $(x_1,x_2)$ of $x$. This means that the electric field $E(x)$ is orthogonal to $B$ and depends only on $(x_1,x_2)$. Assume finally that the initial data $f_0$ only depends on $(x_1,x_2)$ and $(v_1,v_2)$, a property which is therefore propagated by the flow \fref{eq:vlasovJ}.   All these assumptions allow us to restrict ourselves to a $2D\times 2D$ setting and  to rewrite \fref{eq:vlasovJ}  in the form \fref{eq:transport}  with $n=4$, $y=(x,v)\in \R^2\times \R^2$ and
$$
F^\eps(y) = 
\left(
\begin{array}{c}
v \\
\frac{1}{\eps} J v + E(x)
\end{array}
\right)
=
\frac{1}{\eps} G + K \; \mbox{ with } \; 
G(y) = 
\left(
\begin{array}{c}
0 \\
J v 
\end{array}
\right)
\; \mbox{ and } \;
K(y) = 
\left(
\begin{array}{c}
v \\
E(x)
\end{array}
\right).
$$
We now repeat the steps followed for the example of Section \ref{sect:averkin}, starting first with the flow $\Phi_\tau$ (associated with $G$) 
\begin{eqnarray*}
\Phi_\tau(y) = \left(
\begin{array}{c}
x\\
e^{\tau J} v
\end{array}
\right).
\end{eqnarray*}
The time-dependent vector field $K_\tau$ then writes
\begin{eqnarray*}
K_\tau(y) = \left(
\begin{array}{c}
e^{\tau J} v\\
e^{-\tau J} E(x)
\end{array}
\right) = e^{i \tau} \hat K_1(y) + e^{-i \tau} \hat K_{-1}(y) 
\end{eqnarray*}
with 
\begin{eqnarray*}
\hat K_1(y) = \frac12  \left(
\begin{array}{c}
v-i J v\\
E(x)+iJ E(x)
\end{array}
\right)
\quad \mbox{ and } \quad 
\hat K_{-1}(y) =
\frac12  \left(
\begin{array}{c}
v+ iJ v\\
E(x)- iJ E(x)
\end{array}
\right).
\end{eqnarray*}
Formula (\ref{eq:Kseries}) then gives
\begin{eqnarray*}
K^{[1]} &=& \hat K_0 = 0, \\ 
K^{[2]} &=& i [\hat K_1,\hat K_{-1}] = -2 \Im \left((\partial_y \hat K_1) \hat K_{-1} \right) = \left(
\begin{array}{c}
J E\\
\frac12 (\Delta U) J  v
\end{array}
\right)
\end{eqnarray*}
where we used computed successively
\begin{eqnarray*}
\frac{\partial \hat K_1}{\partial y} = 
\frac12 \left(
\begin{array}{cc}
0 & (I-iJ) \\
-(I+iJ) \nabla_x^2 U & 0
\end{array}
\right) 
\end{eqnarray*}
and\footnote{Note that if $S$ is a $2 \times 2$ symmetric matrix then 
$$
J S + S J = 
\left(
\begin{array}{cc}
0 & 1 \\
-1 & 0
\end{array}
\right)
\left(
\begin{array}{cc}
\alpha & \gamma \\
\gamma & \beta
\end{array} 
\right)+
\left(
\begin{array}{cc}
\alpha & \gamma \\
\gamma & \beta
\end{array} 
\right)
\left(
\begin{array}{cc}
0 & 1 \\
-1 & 0
\end{array}
\right)=(\alpha+\beta) J
$$
so that $J \nabla^2 U + \nabla^2 U J = (\Delta U) J$ and 
$$
J S J = 
\left(
\begin{array}{cc}
0 & 1 \\
-1 & 0
\end{array}
\right)
\left(
\begin{array}{cc}
\alpha & \gamma \\
\gamma & \beta
\end{array} 
\right)
\left(
\begin{array}{cc}
0 & 1 \\
-1 & 0
\end{array}
\right)=-\det(S) I.
$$}
\begin{eqnarray*}
\frac{\partial \hat K_1}{\partial y} \; \hat K_{-1}&=&
\frac14
\left(
\begin{array}{c}
(I-iJ)^2 E \\
-(I+iJ) \nabla_x^2 U (I+iJ) v 
\end{array}
\right) \\
&=&
\frac14
\left(
\begin{array}{c}
2(I-iJ) E \\
-\left(\nabla_x^2 U - J \nabla_x^2 U J + i (J \nabla_x^2 U + \nabla_x^2 U J) \right)v 
\end{array}
\right)\\
&=&
\frac14
\left(
\begin{array}{c}
2(I-iJ) E \\
-\left(\nabla_x^2 U + \det( \nabla_x^2 U) I + i \Delta U J \right)v 
\end{array}
\right).
\end{eqnarray*}
At second order in $\eps$,  equation $(i)$ of Corollary \ref{th:main} for $\tilde f(t,\tau,x,v)$  thus has the following form
\begin{eqnarray*}
 \partial_\tau \tilde f + \eps(v - \eps JE) \cdot \nabla_x \tilde f + ((1-\eps^{2} \Delta U) Jv + \eps E) \cdot \nabla_v \tilde f = 0
\end{eqnarray*}
while equation $(ii)$ is simply
\begin{eqnarray*}
 \partial_t \tilde f + \eps J E \cdot \nabla_x \tilde f + \frac{\eps}{2} (\Delta U) J v \cdot \nabla_v \tilde f = 0.
\end{eqnarray*}  
This transport equation  coincides, up to a rescaling in time, with the asymptotic  model derived in \cite{frenod1}.  
We emphasize that, according to Remark \ref{remark-th-S}, these two equations have to be understood  in the approximate sense, which means that they cannot be satisfied
 exactly in general,  but can only be solved approximately allowing errors of order $\eps^2$. 

\subsection{Magnetic field with varying intensity and constant direction} \label{sect:vary}
Over this second phase, we still work in a $2D\times 2D$ setting and keep the same notations as in the previous section. However, we address here the case of a magnetic field with varying intensity $b(x)$ and constant direction $B(x)=(0,0,b(x))^T$.  Note that due to divergence
free property of $B(x)$, the function $b$ depends only on $(x_1,x_2)$.  In order to handle this case of varying intensity $b(x)$, one has to proceed as in Section \ref{sect:varyfreq}. We first immerse
 the problem into an augmented one by
adding a new parametrization variable $s$,  then we derive averaging models at different orders for this augmented problem, and finally eliminate the extra-variable $s$ from the averaged models and show how the original distribution function is recovered.  In order to do so,  we assume that $b(x)$ should not vanish for any $x$ in $\R^2$ and we will make this assumption for the remaining of this section. The augmented distribution function $g(s,t,x,v)$ satisfies 
\begin{eqnarray} \label{eq:vlasovbis}
\pa_s g + \frac{1}{b(x)} \partial_t g + \frac{1}{b(x)} v \cdot \nabla_x g+ 
\left(\frac{1}{\eps} J v -\frac{1}{b(x)} \nabla_x U(x) \right) \cdot \nabla_v g = 0. 
\end{eqnarray}
The original distribution function $f(t,x,v)$ is then viewed as a stationary solution of this evolution equation in $s$.
Denoting $Y=(t,x_1,x_2,v_1,v_2) \in \R^5$ the now {\bf extended} phase-space variable, we equivalently write 
\fref{eq:vlasovbis} as follows
\begin{eqnarray*}
\partial_s g(s,Y) + \check F^\eps(Y)  \cdot \nabla_Y g(s,Y) = 0
\end{eqnarray*}
where 
$$
\check F^\eps(Y) = \left( \begin{array}{c}
\frac{1}{b(x)} \\
\frac{1}{b(x)} v \\
\frac{1}{\eps} J v - \frac{1}{b(x)} \nabla_x U(x)
\end{array}
\right) 
$$
is the  {\bf extended} vector field. We may now resume the derivation of the equations $(i)$ and $(ii)$ of Theorem \ref{th:main}, by first splitting
$\check F^\eps$ into $\check F^\eps = \frac{1}{\eps} \check G + \check K$ with
$$
\check G(Y) = \left( \begin{array}{c}
0 \\
0\\ 
J v
\end{array}
\right)
\quad \mbox{ and } \quad 
\check K(Y) = \frac{1}{b(x)} \left( \begin{array}{c}
1 \\
v \\
- \nabla_x U(x)
\end{array}
\right).
$$
It is clear that $\check G$ now generates a $2 \pi$-periodic flow 
$$
\check \Phi_\tau(Y) = \check \Phi_\tau
\left(\begin{array}{c}
t \\
x \\
v
\end{array}
\right) 
=
\left(\begin{array}{c}
t \\
x \\
e^{\tau J} v
\end{array}
\right)
$$
whose period is independent of the trajectory.
The function $\check K_\tau$ becomes
$$
\check K_\tau (Y) = \frac{1}{b(x)}
\left(\begin{array}{c}
1 \\
e^{\tau J} v \\
e^{-\tau J} E(x)
\end{array}
\right) 
$$
and the corresponding Fourier modes are all vanishing except the modes $1$, $-1$ and $0$ (the additional one w.r.t. the case of a constant field):
\begin{eqnarray*}
\hat K_0(Y) = 
\left(\begin{array}{c}
\frac{1}{b(x)} \\
0 \\
0
\end{array}
\right)\!\!, 
\hat K_1(Y) = 
\frac{1}{2 b(x)}
\left(\begin{array}{c}
0 \\
(I-iJ) v \\
(I+iJ)  E(x)
\end{array}
\right)\!\!,
\hat K_{-1}(Y) = 
\frac{1}{2 b(x)}
\left(\begin{array}{c}
0 \\
(I+iJ) v \\
(I-iJ) E(x)
\end{array}
\right)\!\!.
\end{eqnarray*}
According to Theorem \ref{th:main}, we thus have 
$$
K^{[1]}(Y) = \hat K_0(Y)
$$
and 
\begin{eqnarray*}
K^{[2]} = i \left( [\hat K_1, \hat K_{-1}] + [\hat K_0, \hat K_1-\hat K_{-1}]\right) 
= -2 \Im( [\hat K_0,\hat K_1]) -2 \Im \left((\partial_Y \hat K_1) \hat K_{-1} \right).
\end{eqnarray*}
Omitting the argument $x$ in $E$, $U$ and $b$, and denoting simply $\nabla$ for $\nabla_x$, we have
\begin{eqnarray*}
\frac{\partial \hat K_0}{\partial Y} = 
\frac{-1}{b^2} \left(
\begin{array}{ccc}
0 & \nabla^T b  & 0 \\
0 & 0 & 0 \\
0 & 0 & 0
\end{array}
\right),
\end{eqnarray*}
and 
\begin{eqnarray*}
\frac{\partial \hat K_1}{\partial Y} = 
\frac{1}{2 b^2} \left(
\begin{array}{ccc}
0 & 0 & 0 \\
0 & -(I-iJ)v \, \nabla^T b & b (I-iJ) \\
0 & -b (I+iJ) \nabla^2 U + (I+iJ) \nabla U \,  \nabla^T b & 0
\end{array}
\right) 
\end{eqnarray*}
so that 
\begin{eqnarray*}
(\partial_Y \hat K_1) \hat K_{-1} &=& 
\frac{1}{4 b^3} \left(
\begin{array}{ccc}
0 & 0 & 0 \\
0 & -(I-iJ)v \, \nabla^T b & b (I-iJ) \\
0 & -b (I+iJ) \nabla^2 U + (I+iJ) \nabla U \,  \nabla^T b & 0
\end{array}
\right) 
\left(\begin{array}{c}
0 \\
(I+iJ) v \\
(I-iJ) E
\end{array}
\right) \\
&=& 
\frac{1}{4 b^3}  
\left(\begin{array}{c}
0 \\
-(I-iJ)v \, \nabla^T b \, (I+iJ) v + b (I-iJ)^2 E\\
-b (I+iJ) \nabla^2 U \,  (I+iJ) v + (I+iJ) \nabla U \,  \nabla^T b \, (I+iJ) v
\end{array}
\right) 
\end{eqnarray*}
and finally  
\begin{eqnarray*}
-2 \Im \left((\partial_Y \hat K_1) \hat K_{-1} \right)= \frac{1}{2 b^3}  
\left(\begin{array}{c}
0 \\
   (\nabla b \cdot J v) v 
-   (\nabla b \cdot v) Jv
+2 bJ E\\
 -\eps(\nabla b \cdot v)  J \nabla U - \eps(\nabla b \cdot J v) \nabla U+ b (\Delta U) J v 
\end{array}
\right).
\end{eqnarray*}
Besides, we have 
\begin{eqnarray*}
(\partial_Y \hat K_0) \hat K_1
&=& 
\frac{-1}{2b^3} \left(
\begin{array}{ccc}
0 & \nabla^T b  & 0 \\
0 & 0 & 0 \\
0 & 0 & 0
\end{array}
\right)
\left(\begin{array}{c}
0 \\
(I-iJ) v \\
(I+iJ)  E
\end{array}
\right) \\
&=&
\frac{-1}{2 b^3}
\left(\begin{array}{c}
\nabla b \cdot v - i  \nabla b \cdot J v\\
0 \\
0 
\end{array}
\right)
\end{eqnarray*}
and 
\begin{eqnarray*}
(\partial_Y \hat K_1) \hat K_0
&=& 
\frac{1}{2 b^3} \left(
\begin{array}{ccc}
0 & 0 & 0 \\
0 & -(I-iJ)v \, \nabla^T b & b (I-iJ) \\
0 & -b (I+iJ) \nabla^2 U + (I+iJ) \nabla U \,  \nabla^T b & 0
\end{array}
\right) 
\left(\begin{array}{c}
1 \\
0 \\
0
\end{array}
\right) =0
\end{eqnarray*}
so that 
\begin{eqnarray*}
-2 \Im([\hat K_0,\hat K_1]) = 
\frac{-1}{b^3}
\left(\begin{array}{c}
\nabla b \cdot Jv \\
0 \\
0 
\end{array}
\right).
\end{eqnarray*}
Finally, at first order in $\eps$, we have
\begin{eqnarray*}
\check K^\eps= K^{[1]} + \eps K^{[2]} = 
\frac{1}{b}
\left(\begin{array}{c}
1- \eps  \; \frac{\nabla b \cdot Jv}{b^2} \\
-\eps \; \frac{(\nabla b \cdot v)}{2b^2} J v + \eps \; \frac{(\nabla b \cdot J v)}{2b^2} v - \eps \; \frac{1}{b} J \nabla U \\
-\frac{\eps(\nabla b \cdot v)}{2 b^2}  J \nabla U - \frac{\eps(\nabla b \cdot J v)}{2 b^2} \nabla U+\frac{\eps \Delta U}{2b} J v
\end{array} 
\right) = \left(\begin{array}{c} K_1^\eps\\ K_2^\eps \end{array}\right),
\end{eqnarray*}
and
$$
\check G^\eps= \eps (\check F^\eps -\check K^\eps)= \frac{1}{b}
\left(\begin{array}{c}
0 \\
\eps v \\
 b Jv  - \eps \nabla U
\end{array} 
\right).
$$
Therefore the transport equations on $h$ are at first order in $\eps$:
\begin{eqnarray}
\partial_t h + \frac{\eps}{2b} \left( \frac{\nabla b \cdot J v}{b} v - \frac{\nabla b \cdot v}{b} J v  - 2 J \nabla U\right) \cdot \nabla_x h\nonumber \\
- \frac{\eps}{2b} \left(\frac{\nabla b \cdot v}{b} J \nabla U + \frac{\nabla b \cdot J v}{b} \nabla U-(\Delta U) J v\right) \cdot \nabla_v h = 0,\label{equ-ht-ordre1}
\end{eqnarray}
and
\be\partial_\tau h + \frac{\eps}{b} v\cdot \nabla_x h +  \left( Jv - \frac{\eps}{b} \nabla U\right) \cdot \nabla_v h=0,\label{equ-htau-ordre1}\ee
with the initial condition $h(0,0,y)=f_0(y).$
Similarly the transport equations on $S$ are
\begin{align}
& \partial_t S + \frac{\eps}{2b} \left(\frac{\nabla b \cdot J v}{b} v - \frac{\nabla b \cdot v}{b} J v  - 2 J \nabla U\right) \cdot \nabla_x S \nonumber \\
& - \frac{\eps}{2b} \left(\frac{\nabla b \cdot v}{b} J \nabla U + \frac{\nabla b \cdot J v}{b} \nabla U-(\Delta U) J v\right) \cdot \nabla_v S = b(x)\left( 1+\eps  \; \frac{\nabla b \cdot Jv}{b^2}\right), \label{equ-St-ordre1}
\end{align}
and
\be\partial_\tau S + \frac{\eps}{b} v\cdot \nabla_x S +   \left( Jv - \frac{\eps}{b} \nabla U\right) \cdot \nabla_v S=0, \label{equ-Stau-ordre1}\ee
with the initial condition $S(0,0,y)=0.$  Again, we wish to put the stress on the fact that 
these two truncated models in $h$ and $S$ should be understood
in the sense of Remark \ref{remark-th-S}.

Now we make some important comments on these transport equations. The transport equation \fref{equ-ht-ordre1} coincides with the gyro-kinetic model that has been derived  in \cite{frenod1} in the particular case of constant $b$. It also contains all the terms in the models recently derived in \cite{bostan6} in the case of varying $b=b(x)$ 
when restricted to the  $2D\times 2D$ geometry.  
However, in addition to the fact that our averaged models keep all the variables $(x,v)$, our approach  provides more information through the phase $S$ and the dependence in $\tau$. These informations are necessary to correctly reconstruct  the full  original
distribution function $f$ (and not only the averaged model) at first order in $\eps$. This reconstruction may be performed through the relation $f(t,x,v)= h(t,\tau(t,x,v),x,v)+ O(\eps^2)$  where $\tau(t,x,v)$ is a solution to $\eps\tau = S(t,\tau,x,v)$.
Up to our knowledge, no such construction can be found in the literature.

\subsection{Magnetic field in 3D with varying intensity and varying  direction}

We now consider the transport kinetic equation in its general form (\ref{eq:vlasov}) and in a $3D\times 3D$ setting. This 
means in particular that now we allow  variations of $B$ in both amplitude and direction.
Our aim in this part is to extend our previous approach to this more general case.

We first immerse the model (\ref{eq:vlasov}) into an augmented problem in the unknown
$g (s,t,x,v)$, as follows
\be
\label{eq:vlasov-B-s}
\partial_s g + \frac{1}{|B(x)|}\partial_t g+ \frac{v}{|B(x)|} \cdot \nabla_x g + \left(\frac{E(x)}{|B(x)|} + \frac{1}{\eps} v\times \frac{B(x)}{|B(x)|} \right) \cdot \nabla_v g = 0
\ee
with the initial condition $g(0,t,x,v)= f(t,x,v)$. 
The main interest of this form is that the oscillatory part in the variable $s$ is now driven by the vector field
$v \times\frac{B(x)}{|B(x)|}$, which, as we shall see, generates a periodic flow with a constant period  $2\pi$.  More precisely, the trajectories

$$
\dot x(s) =0, \qquad  \dot v(s)= v(s)\times \frac{B(x(s))}{|B(x(s))|}, \qquad (x(0),v(0))=(x_0,v_0)\in \R^3\times \R^3,
$$
are all periodic with the same period $2\pi$ independently  of $(x_0,v_0).$

In particular the period does not depend on the trajectory although the unit vector $\frac{B(x)}{|B(x)|}$ depends on this trajectory.
Indeed let $e_0$ be a unit vector and let $(e_1,e_2,e_0)$ be  an orthonormal basis such that
$e_0\times e_1=e_2$ and $e_1\times e_2=e_0$.  The matrix  representing the skew-symmetric linear map ${\mathcal J}_{e_0}: v\mapsto v\times e_0$ in the basis $(e_1,e_2,e_0)$,  is simply  ${\mathcal J}=\left(\begin{array}{cc}
J &0\\  0&0\end{array}\right).$  Since $\exp(t{\mathcal J}) $ is $2\pi$-periodic, the flow
$\exp(t{\mathcal J}_{e_0}) $ is $2\pi$-periodic. We now apply our methodology to model \fref{eq:vlasov-B-s}.
Here the vector field $\check F^\eps= \frac{1}{\eps} \check G + \check K$ is given by
$$ \check K(t,x,v) =\frac{1}{|B(x)|}\left( \begin{array}{c} 1 \\ v\\ E(x)
\end{array}\right), \qquad \check G(t,x,v)=\left( \begin{array}{c}  0\\ 0\\ v\times  \frac{B(x)}{|B(x)|}    \end{array}\right).$$
We introduce the following notations
\be
\label{L-P}
e(x)=\frac{B(x)}{|B(x)|}, \quad  {\mathcal J}_{e} v= v\times e, \quad {\mathcal P}_{e} v= (e\cdot v) e, \quad \forall e\in \Sp ^2, v\in \R^3, x\in \R^3.
\ee
Using Theorem \ref{th:series}, the vector field $K_\tau$ can be easily computed to get 
$$
\Phi_\tau (t,x,v)= \left( \begin{array}{c}  t\\ x\\ \exp\left(\tau{\mathcal J}_{e(x)}  \right)  v \end{array}\right).
$$
The following elementary identities
$$ {\mathcal J}_{e}^2 = -I+{\mathcal P}_{e} ,  \qquad {\mathcal J}_{e}{\mathcal P}_{e}= {\mathcal P}_{e}{\mathcal J}_{e}=0
$$
imply that
\begin{align}
\Phi_\tau (t,x,v)& = \left( \begin{array}{c}  t\\ x\\ (\cos\tau) v + (1- \cos\tau){\mathcal P}_{e(x)}v + (\sin\tau){\mathcal J}_{e(x)}   v\end{array}\right)\nonumber\\
& =\left( \begin{array}{c}  t\\ x\\ (\cos\tau) v + (1- \cos\tau)(e(x)\cdot v)e(x) + (\sin\tau)v\times {e(x)}  \end{array}\right) .
\end{align}
We then deduce the expression of the Jacobian matrix $\partial_{t,x,v} \Phi_\tau =\left( \pa_t \Phi_\tau, 
\pa_x \Phi_\tau , \pa_v \Phi_\tau \right)$:
$$
\pa_{t,x,v} \Phi_\tau  = \left( \begin{array}{lll}  1& 0& 0 \\
0 & I &0  \\
0&   R_\tau&
   Q_\tau
\end{array}\right),
$$
where 
$$\begin{array}{ll} R_\tau& = (1- \cos\tau)\pa_x\left( {\mathcal P}_{e(x)} v\right) +(\sin\tau) \pa_x\left( {\mathcal J}_{e(x)} v\right) \\
& = \alpha_0+ \alpha e^{i\tau} + \overline \alpha e^{-i\tau},\\
 Q_\tau&=(\cos\tau) I +(1- \cos\tau){\mathcal P}_{e(x)} + (\sin\tau){\mathcal J}_{e(x)}\\
&=a_0+ a e^{i\tau} + \overline a e^{-i\tau},
\end{array}
$$
and
$$\begin{array}{l} a_0= {\mathcal P}_{e(x)}, \qquad  \alpha_0= \pa_x\left( {\mathcal P}_{e(x)} v\right), \\
2a= I -  {\mathcal P}_{e(x)}-i {\mathcal J}_{e(x)} ,  \\
2\alpha=- \pa_x \left( {\mathcal P}_{e(x)}v+i {\mathcal J}_{e(x)}v  \right). \end{array}  $$
Note that the matrix $R_\tau$ takes care with the so-called curvature terms which are the  terms coming from the space variation of the direction $e(x)$ of the magnetic field.
 In order to compute the inverse of the matrix $\pa_{t,x,v} \Phi_\tau$, we observe
that 
$$
\left(\pa_{t,x,v} \Phi_\tau\right)^{-1}  = \left( \begin{array}{lll}  1& 0& 0 \\
0 & I &0  \\
0&   -Q_\tau^{-1}R_\tau&
Q_\tau^{-1}
\end{array}\right),
$$
which means that we only need to compute $Q_\tau^{-1}$. Using again the identity ${\mathcal J}_{e}^2 = -I+{\mathcal P}_{e}$,  one may check 
$$  
\begin{array}{ll} Q_{\tau}^{-1}& =
(\cos\tau) I +(1- \cos\tau){\mathcal P}_{e(x)} - (\sin\tau){\mathcal J}_{e(x)}, \\
&= a_0 + \overline a e^{i\tau} + a e^{-i\tau}= Q_{-\tau}.
\end{array}
$$
Now we also have 
$$ 
\check K\circ \Phi_\tau (t,x,v)=  \frac{1}{|B(x)|}\left( \begin{array}{c} 1 \\ Q_\tau v \\ E(x)
\end{array}\right), 
$$
and therefore
$$ \check K_\tau (t,x,v)= \frac{1}{|B(x)|} \left( \begin{array}{c}  1 \\
Q_{\tau}v  \\ 
   -Q_{-\tau}R_\tau Q_{\tau}v+
Q_{-\tau}E(x)
\end{array}\right).
$$
One can easily see that the Fourier expansion  of $\check K_\tau$ (in the periodic variable $\tau$)
only contains modes $k\in \Z$ with $|k|\leq 3$.  Note that  we can recover the previous case (in which $B(x)$ had a  constant direction and $(x,v) \in \R^2\times \R^2$) by  taking ${\mathcal P}_{e(x)}v =0$, ${\mathcal J}_{e(x)} \equiv {\mathcal J}= \left(\begin{array}{ccc}J&0\\ 0&0\end{array}\right)$
 and $\alpha=0$, which means that $R_\tau=0$ and $Q_\tau=e^{\tau {\mathcal J} }$.

Although all the Fourier coefficients of $\check K_\tau$ can be derived from this expression, we just give for simplicity the $0^{th}$ mode:
$$\hat K_0(x,v) =  \frac{1}{|B(x)|} \left( \begin{array}{c}  1 \\
{\mathcal P}_{e(x)}v= (e(x)\cdot v) e(x) \\ 
   (\hat K_0)_3
\end{array}\right) = K^{[1]},$$
with
\begin{align*} 
(\hat K_0)_3(x,v) &=  a_0 E(x) - ( a_0\alpha_0a_0+ a_0\alpha\overline a+ a_0\overline \alpha a +
\overline a \alpha_0 \overline a + \overline a \overline \alpha a_0 + a \alpha_0 a + a \alpha a_0) v \\
&=   {\mathcal P}_{e(x)} E(x) - 
\left[  4 {\mathcal P}_{e(x)} \pa_x \left( {\mathcal P}_{e(x)}v\right) {\mathcal P}_{e(x)} + \frac{1}{2}  {\mathcal P}_{e(x)} \pa_x \left( {\mathcal J}_{e(x)}v\right) {\mathcal J}_{e(x)}\right. \\ & \left.
- \frac{1}{2} {\mathcal J}_{e(x)} \pa_x \left( {\mathcal P}_{e(x)}v\right) {\mathcal J}_{e(x)}  - \frac{1}{2}  {\mathcal J}_{e(x)} \pa_x \left( {\mathcal J}_{e(x)}v\right) {\mathcal P}_{e(x)} \right. \\
&\left. 
- {\mathcal P}_{e(x)} \pa_x \left( {\mathcal P}_{e(x)}v\right)- 
 \pa_x \left( {\mathcal P}_{e(x)}v\right){\mathcal P}_{e(x)} + \frac{1}{2}  \pa_x \left( {\mathcal P}_{e(x)}v\right)
\right] v.
\end{align*}
We then deduce the vector field $G^\eps$ at the $0^{th}$ order in $\eps$:
$$\begin{array}{ll} G^{[1]} = \eps (\check F^\eps -K^{[1]}) + {\cal O}(\eps)&= \frac{1}{|B(x)|}
\left(\begin{array}{c}
0 \\
\eps \left (v -{\mathcal P}_{e(x)}v\right)  \\
 |B(x)| {\mathcal L}_{e(x)}v + \eps \left( E(x)- (\hat K_0)_3(x,v) \right)
\end{array} 
\right) + {\cal O}(\eps)\\
&= 
\left(\begin{array}{c}
0 \\
0 \\
{\mathcal L}_{e(x)}v 
\end{array} 
\right) + {\cal O}(\eps).
\end{array}
$$

The averaged model  at the $0^{th}$ order in $\eps$ can now be written in terms of   $h(t,\tau,x,v)$ and
$S(t,\tau,x,v)$. We have 
$$  
\begin{array}{l}
\pa_t h + \left( \frac{B(x)}{|B(x)|}\cdot v\right) \frac{B(x)}{|B(x)|} \cdot \nabla_x h + (\hat K_0)_3(x,v)\cdot \nabla_v h =0,\\
\pa_\tau h  + \left(v\times \frac{B(x)}{|B(x)|}\right) \cdot \nabla_v h =0,
\end{array}
$$
and 
$$  
\begin{array}{l}
\pa_t S + \left( \frac{B(x)}{|B(x)|}\cdot v\right) \frac{B(x)}{|B(x)|}\cdot \nabla_x S + (\hat K_0)_3(x,v)\cdot \nabla_v S =|B(x)|,\\
\pa_\tau S  + \left(v\times \frac{B(x)}{|B(x)|}\right) \cdot \nabla_v S =0,
\end{array}
$$
with the initial conditions: $ h(0,0,x,v)= f_0(x,v)$ and $S(0,0,x,v)=0.$
Note that in the particular case where $B(x)$ has a constant direction $B(x)=b(x)e_0= (0,0,b(x))^T$, we get
$$  
\begin{array}{l}
\pa_t h + \ v_{\parallel} \pa_{x_{\parallel}} h + E_{\parallel} \pa_{v_{\parallel}} h =0,\\
\pa_\tau h  + Jv_{\perp} \cdot \pa_{v_{\perp}} h =0.
\end{array}
$$
and 
$$  
\begin{array}{l}
\pa_t S + \ v_{\parallel} \pa_{x_{\parallel}} S + E_{\parallel} \pa_{v_{\parallel}} S =b(x),\\
\pa_\tau S  + Jv_{\perp} \cdot \pa_{v_{\perp}} S =0,
\end{array}
$$
where  we used the standard notations $v_{\parallel}= v\cdot e_0,\  E_{\parallel}= E\cdot e_0$ and 
$$ v= (v_1,v_2,v_\parallel)= (v_{\perp},v_{\parallel}), \quad E=(E_1,E_2,E_\parallel)=(E_{\perp},E_{\parallel}),\quad \pa_{v_{\perp}} h =  (\pa_{v_1}h, \pa_{v_2}h), $$
and the same notations for the space variable $x$. Observe that the exact solution of the 
two equations for $S$ (for the $0^{th}$ order in $\eps$) is simply $S(t,\tau,x,v)= b(x) t$.

The averaged equations at the first order in $\eps$ can also be derived  in the case of a magnetic field $B(x)$ with constant direction $B(x)= (0,0,b(x))^T$, with $b(x)>0$.
In this case we have $R_\tau=0$, $e(x)$ is the constant unit vector $e_0$, $|B(x)|= b(x),$ and therefore
$$ 
 \check K_\tau (t,x,v)=\frac{1}{b(x)} \left( \begin{array}{c}  1 \\
Q_{\tau}v  \\ 
Q_{-\tau}E(x)
\end{array}\right).
$$
The non-zero Fourier modes  in $\tau$ of this quantity $\check K_\tau$  are
$$ 
\hat K_0= \frac{1}{b} \left( \begin{array}{l}  1 \\
a_0 v  \\ 
 a_0 E
\end{array}\right), \quad \hat K_1 =\frac{1}{b} \left( \begin{array}{l}  0 \\
a v  \\ 
\overline a E
\end{array}\right), \quad 
 \hat K_{-1} =\frac{1}{b} \left( \begin{array}{l}  0 \\
\overline a v  \\ 
 a E
\end{array}\right).
$$
The computation of  $K^\eps$ at first order in $\eps$ can then be derived from
 Theorem \ref{th:series} as follows.
 We know from Theorem \ref{th:series} that $ \check K^\eps = K^{[1]}+ \eps K^{[2]} $ with
 $$
 K^{[1]}= \hat K_0, \quad K^{[2]} = -2\Im \left( (\pa_Y \hat K_1 ) \hat K_{-1}\right) - 
 2\Im \left([ \hat K_0,\hat K_1] \right).
 $$
 Since
 $$
 \pa_Y \hat K_1= \left( \begin{array}{ccc}  0& 0& 0 \\
0 & -a \left(v\otimes \frac{\nabla b}{b^2}\right) &\frac{a}{b}  \\
0&   \overline a \pa_x \left( \frac{E}{b}\right) &
0
\end{array}\right), \quad \pa_Y \hat K_0=\left( \begin{array}{ccc}  0& -\frac{(\nabla b)^T}{b^2}& 0 \\
0 & -a_0 \left(v\otimes \frac{\nabla b}{b^2}\right) &\frac{a_0}{b}  \\
0&   a_0 \pa_x \left( \frac{E}{b}\right) &
0
\end{array}\right)
 $$
we get
$$ 2\Im \left( (\pa_Y \hat K_1 ) \hat K_{-1}\right)= \frac{1}{2b} 
\left( \begin{array}{c}  0\\  -(I-\mathcal P) \left( v\otimes \frac{\nabla b}{b^2}\right) {\mathcal J}v 
+ {\mathcal J} \left( v\otimes \frac{\nabla b}{b^2}\right)(I-\mathcal P) v
  -\frac{2}{b}{\mathcal J} E\\
(I-\mathcal P)  \pa_x\left(\frac{E}{b}\right){\mathcal J} v +  {\mathcal J}\pa_x\left(\frac{E}{b}\right)(I-\mathcal P) v  
\end{array}\right)
$$
where we have denoted  $\mathcal P = \mathcal P_{e_0}$ and  $\mathcal J = \mathcal J_{e_0}$. We also have
$$ 2\Im \left([ \hat K_0,\hat K_1] \right)= \frac{1}{b} 
\left( \begin{array}{c} {\mathcal J}v \cdot  \frac{\nabla b}{b^2} \\  \mathcal P \left( v\otimes \frac{\nabla b}{b^2}\right) {\mathcal J}v 
- {\mathcal J} \left( v\otimes \frac{\nabla b}{b^2}\right)\mathcal P v
  \\
-\mathcal P \pa_x\left(\frac{E}{b}\right){\mathcal J} v -  {\mathcal J}\pa_x\left(\frac{E}{b}\right)\mathcal P v  
\end{array}\right),
$$
therefore
$$ K^{[2]}= \frac{1}{b} 
\left( \begin{array}{c} - {\mathcal J}v \cdot  \frac{\nabla b}{b^2}  \\ 
 \frac{1}{2} \left({\mathcal J} v\cdot \frac{\nabla b}{b^2}\right) (I-3\mathcal P) v 
- \frac{1}{2} \left( (I-3\mathcal P)v\cdot \frac{\nabla b}{b^2}\right){\mathcal J} v
  +\frac{1}{b}{\mathcal J} E\\
- \frac{1}{2}(I-3\mathcal P)  \pa_x\left(\frac{E}{b}\right){\mathcal J} v -  \frac{1}{2}{\mathcal J}\pa_x\left(\frac{E}{b}\right)(I-3\mathcal P) v  
\end{array}\right),
$$
and
\begin{align*}
\check K^\eps &= K^{[1]}+ \eps K^{[2]}+ O(\eps^2) \\
   & =\frac{1}{b}\left( \begin{array}{c}  1- \eps {\mathcal J}v \cdot  \frac{\nabla b}{b^2}  \\ 
 v_{\parallel} e_0 +  \frac{\eps}{2}\left[ \left({\mathcal J} v\cdot \frac{\nabla b}{b^2}\right) (I-3\mathcal P) v 
-  \left( (I-3\mathcal P)v\cdot \frac{\nabla b}{b^2}\right){\mathcal J} v\right]
  +\frac{\eps}{b}{\mathcal J} E\\
 E_{\parallel} e_0- \frac{\eps}{2}\left[(I-3\mathcal P)  \pa_x\left(\frac{E}{b}\right){\mathcal J} v + {\mathcal J}\pa_x\left(\frac{E}{b}\right)(I-3\mathcal P) v \right] 
\end{array}\right) +  O(\eps^2).
\end{align*}
We finally deduce the field $G^\eps$ at first order in $\eps$
$$
\check G^\eps = \eps (\check F^\eps - K^{[1]}- \eps K^{[2]}) + O(\eps^2) = \left(\begin{array}{c} 0\\0\\ \mathcal J v\end{array}\right) +
\frac{\eps}{b} \left(\begin{array}{c} 0\\v_{\perp }\\ E_{\perp}\end{array}\right) + O(\eps^2).
$$
 Therefore, the evolution in time $t$ of   $h$ at the first order in $\eps$ is driven by the following equation (with the above described notations) 

$$  
\begin{array}{l}
\ds \left[1-\eps  J v_\perp  \cdot \frac{\pa_{x_\perp} b}{b^2} \right] \ds \pa_t h    +  \ v_{\parallel}\left[ 1-\eps  J v_\perp  \cdot \frac{\pa_{x_\perp} b}{b^2} \right] \pa_{x_\parallel} h 
  + \left[ E_{\parallel} +   \eps \pa_{x_\perp}\left(\frac{E_\parallel}{b}\right) \cdot v_\perp \right] \pa_{v_\parallel} h  \\
  \hspace{2.1cm}\ds 
 - \frac{\eps}{2b} \left [ |v_\perp|^2 \frac{\pa_{x_\perp} b}{b} - 2 JE_{\perp}\right]  \cdot \pa_{x_{\perp}} h\\ 
\hspace{2.1cm} \ds +\frac{\eps}{2}  \left[  \left(\frac{\pa_{x_\perp}b}{b^2} \cdot JE_\perp\right) v_\perp +2v_\parallel \pa_{x_\parallel}\left(\frac{E_\perp}{b}\right) - 
  \pa_{x_\parallel}\left(\frac{E_\perp}{b}\right)  Jv_\perp \right]\cdot \pa_{v_\perp} h =0,
\end{array}
$$
which simplifies into
\be 
\label{eq-h-t-order1}
\begin{array}{l}
\ds \pa_t h    +  \ v_{\parallel} \pa_{x_\parallel} h 
 +  \left[ E_{\parallel} + \eps E_{\parallel} J v_\perp  \cdot \frac{\pa_{x_\perp} b}{b^2} + \eps \pa_{x_\perp}\left(\frac{E_\parallel}{b}\right) \cdot v_\perp \right] \pa_{v_\parallel} h \\ 
 \hspace{2.1cm}\ds - \frac{\eps}{2b} \left [ |v_\perp|^2 \frac{\pa_{x_\perp} b}{b} - 2 JE_{\perp}\right]  \cdot \pa_{x_{\perp}} h\\ 
\hspace{2.1cm}\ds +\frac{\eps}{2}  \left[  \left(\frac{\pa_{x_\perp}b}{b^2} \cdot JE_\perp\right) v_\perp +2v_\parallel \pa_{x_\parallel}\left(\frac{E_\perp}{b}\right) - 
  \pa_{x_\parallel}\left(\frac{E_\perp}{b}\right)  Jv_\perp \right]\cdot \pa_{v_\perp} h =0.
\end{array}
\ee
Note that  we have used the identity $\nabla_x \cdot B=0$ which implies that $b(x)= b(x_\perp)$.   
This provides  an asymptotic model which is identical to the one recently derived  in \cite{bostan6} or, up to a rescaling in time,  to the one derived in \cite{degond2}. However our approach provides more informations since this equation still contains all the original variables $(x,v)$ of the distribution function and has to be coupled with an equation describing its dependence on a periodic variable $\tau$ which has to fit  with a suitable phase function $S$. As we shall see, this equation  in $\tau$ will provide a suitable initial data for equation \fref{eq-h-t-order1}.
The second equation on $h$ writes
\be
\label{eq-h-tau-order1} \pa_\tau h+ Jv_\perp\cdot \pa_{v_\perp} h+  \frac{\eps}{b}  v_\perp \cdot \pa_{x_\perp} h  + \eps \frac{E_\perp}{b} \cdot
\pa_{v_\perp}h=0.
\ee
The system of the two equations (\ref{eq-h-t-order1}-\ref{eq-h-tau-order1})   is subjected to the initial data $h(0,0,x,v)=f_0(x,v).$

Once again, we recall that system (\ref{eq-h-t-order1}-\ref{eq-h-tau-order1})  with  initial condition $h(0,0,x,v)=f_0(x,v)$ is only valid up to  $\eps^2$ terms, and  solutions to this system have to be understood in the sense of Remark \ref{remark-th-S}. 

Similarly the equations on $S$ are

\be 
\label{eq-S-t-order1} 
\begin{array}{l}
\ds \pa_t S    +  \ v_{\parallel} \pa_{x_\parallel} S 
 +  \left[ E_{\parallel} + \eps E_{\parallel} J v_\perp  \cdot \frac{\pa_{x_\perp} b}{b^2} + \eps \pa_{x_\perp}\left(\frac{E_\parallel}{b}\right) \cdot v_\perp \right] \pa_{v_\parallel} S \\ 
 \hspace{2.5cm}\ds - \frac{\eps}{2b} \left [ |v_\perp|^2 \frac{\pa_{x_\perp} b}{b} - 2 JE_{\perp}\right]  \cdot \pa_{x_{\perp}} S\\
  \hspace{2.5cm}\ds+\frac{\eps}{2}  \left[  \left(\frac{\pa_{x_\perp}b}{b^2} \cdot JE_\perp\right) v_\perp +2v_\parallel \pa_{x_\parallel}\left(\frac{E_\perp}{b}\right) - 
  \pa_{x_\parallel}\left(\frac{E_\perp}{b}\right)  Jv_\perp \right]\cdot \pa_{v_\perp} S\\
   \hspace{2.5cm}\ds  = b+\eps  J v_\perp  \cdot \frac{\pa_{x_\perp} b}{b},
\end{array}
\ee
and 
\be 
\label{eq-S-tau-order1}  \pa_\tau S + Jv_\perp \cdot \pa_{v_\perp} S+  \frac{\eps}{b}  v_\perp \cdot \pa_{x_\perp} S  + \eps \frac{E_\perp}{b} \cdot
\pa_{v_\perp}S=0,\ee
with the initial data $S(0,0,x,v)=0.$

We now observe that $S(t,\tau,x,v)= b(x) t + O(\eps)$, and therefore it is more convenient to write these equations in terms
of $$ \tilde S(t,\tau,x,v)=\frac{S(t,\tau,x,v)-b(x)t}{\eps}$$
and get
$$  
\begin{array}{l}
\ds \pa_t \tilde S    +  \ v_{\parallel} \pa_{x_\parallel} \tilde S 
 +  \left[ E_{\parallel} + \eps E_{\parallel} J v_\perp  \cdot \frac{\pa_{x_\perp} b}{b^2} + \eps \pa_{x_\perp}\left(\frac{E_\parallel}{b}\right) \cdot v_\perp \right] \pa_{v_\parallel} \tilde S \\ 
 \hspace{2.5cm}\ds - \frac{\eps}{2b} \left [ |v_\perp|^2 \frac{\pa_{x_\perp} b}{b} - 2 JE_{\perp}\right]  \cdot \pa_{x_{\perp}} \tilde S\\
 \hspace{2.5cm}\ds+\frac{\eps}{2}  \left[  \left(\frac{\pa_{x_\perp}b}{b^2} \cdot JE_\perp\right) v_\perp +2v_\parallel \pa_{x_\parallel}\left(\frac{E_\perp}{b}\right) - 
  \pa_{x_\parallel}\left(\frac{E_\perp}{b}\right)  Jv_\perp \right]\cdot \pa_{v_\perp} \tilde S\\
   \hspace{2.5cm}\ds =   J v_\perp  \cdot \frac{\pa_{x_\perp} b}{b},
\end{array}
$$
and 
$$  \pa_\tau \tilde S+ Jv_\perp\cdot \pa_{v_\perp} \tilde S+\eps \frac{t}{b}  v_\perp \cdot \pa_{x_\perp} b +  \frac{\eps}{b}  v_\perp \cdot \pa_{x_\perp} \tilde S  + \eps \frac{E_\perp}{b} \cdot
\pa_{v_\perp}\tilde S=0,$$
with the initial data   $\tilde S(0,0,x,v)=0.$
We then recover the solution $f$ by the relation
$$
f(t,x,v) = h(t,\tau(t,x,v),x,v),
$$
where $(t,x,v) \mapsto \tau(t,x,v) \in \R$ is implicitly defined (locally) from $\tilde S$ by the equation
$$
 \tau(t,x,v)= \frac{b(x)t}{\eps}+ \tilde S(t,\tau(t,x,v),x,v).
$$

%

%

\end{document}